\documentclass[3p, number,sort&compress,times,fleqn]{elsarticle}      

\usepackage{./csml}

\usepackage{braket}

\usepackage[colorlinks = true,
linkcolor = blue,
urlcolor  = blue,
citecolor = blue,
anchorcolor = blue]{hyperref}

\usepackage{tikz}
\usetikzlibrary{quantikz2}
\usetikzlibrary{calc,positioning,arrows.meta,fpu}
\tikzset{operator/.append style={fill opacity=0.85}}

\usepackage{nicematrix}

\usepackage{placeins}

\begin{document}

\begin{frontmatter}

\title{IQFEM: A quantum finite element method for heterogeneous problems with immersed boundaries}

\author{Shehara Perera}
\author{Yiren Wang}
\author{Fehmi Cirak\corref{cor1}}
\ead{f.cirak@eng.cam.ac.uk}

\cortext[cor1]{Corresponding author}

\address{Department of Engineering, University of Cambridge, Cambridge, CB2 1PZ, UK }

\begin{abstract}
Quantum computing offers the potential to solve computational mechanics problems with a more favourable complexity scaling than classically possible. To realise this potential, a quantum-compatible formulation of the finite element method is essential. We present such a formulation for Poisson problems with spatially varying coefficients on general domains discretised by an immersed uniform Cartesian grid. The resulting banded sparse system matrices are encoded as the sum of diagonal matrices with polynomial entries and their products with shift matrices. The sum of the matrices is formed using the linear combination of unitaries (LCU) technique. The resulting block-encoded system matrix can be implemented using a polylogarithmic number of elementary gates. For immersed boundaries, the block-encoded system matrix is further processed using a diagonal indicator matrix with zeros and ones. The indicator matrix is defined in terms of an integer-valued level-set function composed of simple geometric primitives via Boolean set operations. The required integer operations are computed on the fly using quantum arithmetic. The resulting block-encoded linear system of equations is solved using quantum singular value transformation (QSVT). The overall gate count, including the QSVT, is~$O(N_{\text{tot}}^{2/d} \operatorname{polylog}(N_{\text{tot}}))$, where~$d$ is the problem dimension and~$N_{\text{tot}}$ the number of grid points. In numerical examples that confirm the theoretical complexity estimates, the solutions converge approximately linearly for immersed problems and quadratically for simply connected domains. 
\end{abstract}


\end{frontmatter}

%
\section{Introduction \label{sec:intro}}
\begin{figure}[!htbp]
    \centering
		\includegraphics{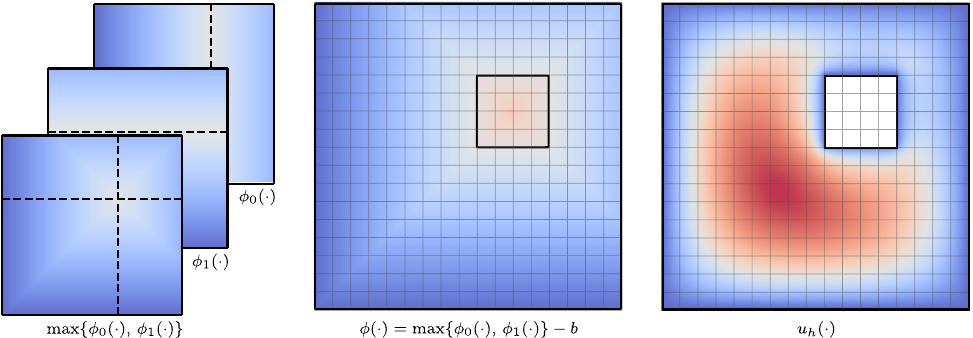} 
        \caption{{\em IQFEM} analysis of a multiply connected domain with an internal hole.  The geometry of the hole, with half-width~$b$, is described by a level-set function~$\phi$ composed of two simpler level-set functions~$\phi_0$ and~$\phi_1$ via constructive solid geometry.  In {\em IQFEM}, the integer-valued level-set functions are evaluated on the fly using quantum arithmetic.  The linear system of finite element equations is solved with quantum singular value transformation, yielding the discrete solution~$u_h$.}\label{fig:intro}
\end{figure}

%
Quantum algorithms for solving linear systems of finite element equations scale, according to theoretical estimates, polylogarithmically in the system dimension and polynomially in the condition number~\cite{harrow2009,childs2017,gilyen2019quantum}. The prospect of achieving an exponential speedup over classical algorithms with polynomial complexity represents a paradigm shift. However, the complexity of quantum linear systems (QLS) algorithms depends on crucial details of data encoding and solution readout~\cite{aaronson2015read,montanaro2016quantum,morales2024quantum}. In this work, we focus on efficient block encoding of finite element matrices, a key bottleneck to achieving the potential speedup. Theoretical complexity estimates for QLS algorithms are typically given in terms of the number of oracle queries, namely, a black-box unitary that provides access to the system matrix. The complexity estimates, therefore, tacitly assume that the oracle itself can be implemented with only polylogarithmically many elementary gates in the system size. This is usually difficult to achieve for a system matrix resulting from a fully unstructured finite element discretisation. We therefore consider finite element discretisations on structured Cartesian grids and exploit their regularity in constructing the required circuits. Specifically, we consider domains of arbitrary shape using the immersed finite element approach, also known as an embedded or cut finite element approach~\cite{febrianto2024three,main2018shifted,Burman2015,Sanches:2011aa,Ruberg:2010aa,Parvizian2007}.  We refer to the resulting immersed quantum finite element method as~{\em IQFEM} (Figure~\ref{fig:intro}).

The efficient encoding of finite element system matrices has motivated a range of approaches. Quantum spectral methods reformulate the problem in a basis in which the system matrix has a simple diagonal or block-diagonal structure~\cite{childs2020quantum,liu2024towards,lubasch2025quantum,febrianto2025quantum,wang2026qafe,gimenez2026quantum}. Such a reformulation is generally not possible for problems with heterogeneous coefficients and domains of arbitrary shape. Therefore, in {\em IQFEM} we directly block-encode the system matrices arising from immersed finite element discretisations on structured Cartesian grids. Finite element matrices are usually non-unitary and must therefore be block-encoded as a sub-block of a larger, suitably chosen unitary matrix. Such block encodings are commonly constructed as a sequence of unitary transformations that combine two-dimensional rotations on selected subspaces of the quantum Hilbert space.  The rotation angles depend on the entries of the system matrix. Efficiency requires that the number of gates to select the subspaces, i.e., $\mathit{CNOT}$ gates, and to apply the rotations, i.e.~$R_Y$ gates, depends polylogarithmically on the matrix dimension. Approaches that make no assumptions about the sparsity pattern or the entries of the matrix usually lead to quadratic gate complexity in system size~$N_{\text{tot}}$ for a matrix of dimension $N_{\text{tot}} \times N_{\text{tot}}$~\cite{bergholm2005quantum,shende2005synthesis,arora2025implementation}. For sparse matrices, this cost can be reduced via circuit synthesis techniques that combine~$R_Y$ gates and neglect those with small angles~\cite{mottonen2004quantum,camps2022fable}. This reduction alone does not, in general, achieve polylogarithmic scaling.

The finite element system matrices for problems discretised using structured grids have a particularly simple banded sparsity structure that can be exploited for their quantum encoding. Furthermore, the matrices have repeated entries when the elements are of the same size, and the problem has constant coefficients. For instance, a periodic one-dimensional Poisson problem discretised with linear basis functions has a tridiagonal circulant system matrix with entries proportional to $-1$, $2$, and $-1$. The circulant system matrix can be represented as a weighted linear combination of an easily block-encoded diagonal matrix and two periodic shift matrices~\cite{camps2024explicit}. The shift matrices represent the off-diagonal terms and are unitary permutations. The weighted sum of the block-encoded diagonal matrix and the unitary shift matrices can be obtained using the linear combination of unitaries (LCU) technique~\cite{childs2012hamiltonian,herbert2026quantum}. The unitary for the entire quantum circuit serves as a block-encoding of the system matrix. Extensions and variations of the sketched approach are considered in~\cite{sunderhauf2024block,kharazi2025explicit}. 

In this paper, we develop quantum algorithms with polylogarithmic gate complexity for the block encoding of system matrices arising from problems with heterogeneous coefficients and domains of general shape.  The specific problems considered are the one- and two-dimensional heterogeneous Poisson problems discretised on immersed Cartesian grids using linear basis functions. The formulation can be readily extended to three- and higher-dimensional problems. The spatial distribution of the diffusivity coefficient is assumed to be polynomial and is approximated as constant within each element. Consequently, the entries of the system matrix become polynomials of the row indices, which correspond to the coordinates of the grid points. The resulting polynomial dependence of the matrix entries on the grid indices enables the extension of the approach of Camps et al.~\cite{camps2024explicit} from constant-coefficient to heterogeneous problems. Specifically, in one dimension, the tridiagonal stiffness matrix is the sum of a diagonal and two off-diagonal matrices with polynomial entries. The two off-diagonal matrices are products of diagonal matrices with two shift operators.  We encode the polynomial entries of the three diagonal matrices using an approach with polylogarithmic complexity~\cite{woerner2019quantum,rosenkranz2025quantum}. The sum of the block-encoded diagonal and off-diagonal matrices is determined using the LCU technique. For instance, in two dimensions, we decompose the system matrix into one diagonal matrix and four off-diagonal matrices with bivariate polynomial entries. The off-diagonal matrices are products of diagonal and shift matrices. We incorporate non-periodic boundary conditions by altering the shift operators.

We block-encode the system matrices of domains embedded in Cartesian grids by extending our approach for heterogeneous problems. Without loss of generality, we focus on multiconnected domains with an internal hole subject to either homogeneous Dirichlet or Neumann boundary conditions. To describe the problem domain, we derive a block-encoded diagonal masking matrix with entries equal to zero outside the domain and one within it. The system matrix for homogeneous Dirichlet boundary conditions on the hole is obtained by multiplying the block-encoded system matrix of the heterogeneous problem from the left and right by the block-encoded diagonal masking matrix. The resulting approach has linear convergence in the~$L^2$-norm~\cite{ramiere2008convergence,main2018shifted}. The system matrix is further modified to impose homogeneous Neumann boundary conditions. Specifically, the diagonal entries of the system matrix corresponding to the boundary grid points are scaled according to the number of adjacent cells outside the domain. This approach can be motivated by a finite volume treatment of the Neumann boundaries~\cite{eymard2000finite}.  Although the obtained system matrices for multiconnected domains have zero rows and columns, this does not affect the final solution, as the corresponding singular values are suppressed during the quantum solution stage. 

The entries of the masking matrix required to block-encode the system matrices of multiconnected domains are given by an indicator function defined in terms of a level-set function (Figure~\ref{fig:intro}). Efficient encoding of the level-set function with a polylogarithmic number of gates is key to the proposed approach's efficiency. A polynomial approximation of the level-set function generally requires a high polynomial degree and, consequently,  a large number of gates. Instead, we compute the level-set function on the fly using a constructive solid geometry approach~\cite{Ricci1973}, in which the domain is recursively composed of simple geometric primitives, such as half-spaces. The level-set and indicator functions are represented using binary encoding and computed using quantum arithmetic.  
In binary encoding, an integer or floating-point number is represented using computational basis states; for example, the integers $0, \, \dotsc, \, 15$  are represented by four qubits, i.e., $0000, \, \dotsc, \,1111$,  analogous to a four-bit representation in classical computing. We use quantum Fourier transform (QFT)-based quantum arithmetic~\cite{draper2000addition,ruiz2017quantum,wang2025comprehensive} to compute the level-set and indicator functions with a polylogarithmic number of gates. In our present implementation, the embedded domain boundary is approximated by element edges, yielding a staircase approximation of the boundary, and no cut elements are considered. This approximation can potentially be improved by choosing a finer subgrid for capturing the contributions of cut elements~\cite{febrianto2024three}.  

We solve the block-encoded linear system of equations using the quantum singular value transformation (QSVT) approach. QSVT is a quantum approach that represents matrix polynomials as products of block-encoded matrices interleaved with parametrised signal-processing operators~\cite{gilyen2019quantum,martyn2021grand}.  Each signal-processing operator depends on a single parameter, i.e., a phase factor, and the set of phase factors depends on the coefficients of the polynomial under consideration. The phase factors are not unique, but efficient algorithms exist for determining them~\cite{dong2021qsppack}.  In approximating the matrix inverse using QSVT, we first approximate the function~$x^{-1}$ by a polynomial while suppressing the singularity near the origin and then compute the corresponding phase factors.  The same phase factors are used to approximate the inverse of the block-encoded system matrix. The query complexity, or the number of times the block-encoded matrix is evaluated in QSVT when approximating the matrix inverse, is given as $O(\kappa \log (\kappa/\epsilon))$, where $\kappa$ is the condition number and $\epsilon$ the approximation error~\cite{gilyen2019quantum}. Assuming that~$\kappa \sim h^{-2} \sim N_{\text{tot}}^{2/d}$, where $h$ is the element size, $N_{\text{tot}}$ the total number of grid points and $d$ the space dimension, the query complexity is $O(N_{\text{tot}}^{2/d} \log (N_{\text{tot}}^{2/d}/\epsilon))$. This scaling motivates the preconditioning of the system matrix, especially for $d<3$, which is not considered in the present paper. 

The remainder of the paper is organised as follows. In Section~\ref{sec:problem}, we introduce the finite element system matrices for one- and two-dimensional Poisson problems discretised on structured Cartesian meshes using linear basis functions. In Section~\ref{sec:sysmat}, we discuss the quantum encoding of the system matrices of problems with a heterogeneous diffusion coefficient. Section~\ref{sec:immersed} extends the formulation to multiconnected domains with an internal hole embedded in Cartesian meshes, with particular attention to the computation of level-set and indicator functions and to the treatment of Dirichlet and Neumann boundary conditions. In Section~\ref{sec:lsolve}, we describe the solution of the resulting linear systems using QSVT. Numerical examples assessing accuracy and scaling are presented in Section~\ref{sec:examples}. The paper concludes with a discussion of implications, limitations, and directions for future work in~{\em IQFEM}.  Three appendices provide additional details.

%
\section{Model problem: heterogeneous diffusion \label{sec:problem}}
%
In this section, we recall the finite element discretisation of the Poisson equation with a non-homogeneous diffusion coefficient and give the explicit expressions for the system matrices in one and two dimensions. The domain is partitioned using a structured mesh with simplicial elements, and the weak form of the Poisson equation is discretised with standard linear basis functions. 
%
\subsection{Governing equation and discretisation \label{sec:govAndDisc}}
%
On the  domain~\mbox{$\Omega \subseteq (0, \, L)^d $}, with $d \in \{ 1, \, 2 \}$, which may contain holes, we consider the governing equation 
\begin{subequations}	 \label{eq:poisson}
\begin{alignat}{2}
	- \nabla \cdot  \left ( \mu  \nabla u\right ) &=  f   \qquad &&\text{in }  \Omega  \, ,  \\ 
	u  &= 0 &&\text{on } \Gamma_D \, ,   \\ 
	\mu \nabla u  \cdot n  &=  t &&\text{on } \Gamma_N \, , 
\end{alignat}
\end{subequations}
where~$u$ is the unknown solution field,~$\mu \in L^\infty (\Omega)$ is the diffusion coefficient and~$ f \in L^2 (\Omega)$ is the source term.  The spaces~$ L^\infty (\Omega)$ and~$ L^2 (\Omega)$ denote the standard spaces of essentially bounded functions and square-integrable functions, respectively.  The prescribed flux on the Neumann boundary~$\Gamma_N$  is~$t$, and the prescribed solution on the Dirichlet boundary~$\Gamma_D$ is zero.  Without loss of generality, the external boundary is assumed to be entirely Dirichlet, and any hole boundaries are either Dirichlet or Neumann.

The weak form of the governing equation~\eqref{eq:poisson} is stated as: find~$u \in H_0^1 (\Omega)$ such that 
\begin{equation}  \label{eq:weak_form}
	\int_\Omega \mu \nabla u \cdot  \nabla v  \D \Omega  = \int_\Omega f v \D \Omega   + \int_{\Gamma_N} t v  \D{\Gamma}  \, , \qquad \forall v \in H_0^1 (\Omega)  \, ,
\end{equation}
where $H_0^1 (\Omega) = \{ v \in  H^1(\Omega) \vert v = 0 \text{ on }  \Gamma_D \}$, and~$H^1(\Omega)$ is the Sobolev space of square-integrable functions whose first derivative is also square-integrable. 
 
\begin{figure}[!htb]
    \centering
    	\subfloat[][One dimension \label{fig:1d-mu-grid}] {
		\includegraphics[scale=0.92, trim={4 0 2 0}, clip]{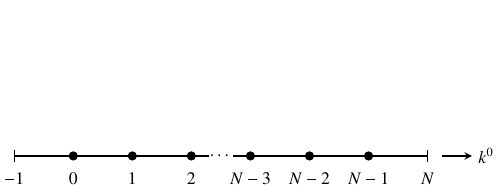} 
	}
   	\hfill
	\subfloat[][Two dimensions  \label{fig:2d-mu-grid}] {
		\includegraphics[scale=0.92, trim={3 0 7 0}, clip]{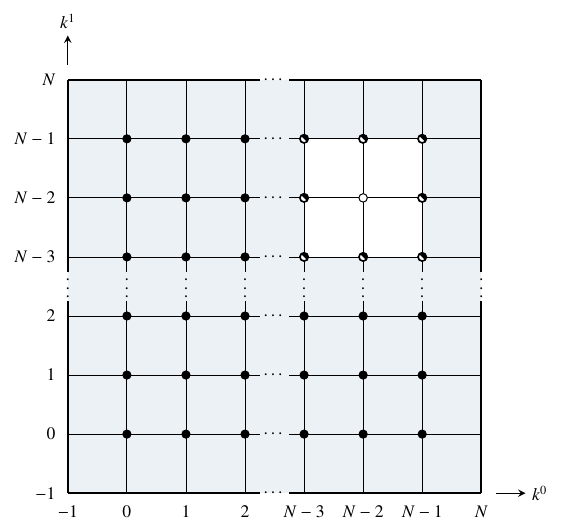} 
        }
    \caption{Discretisation of the domain~$\Omega \subseteq (0, \, L)^d$, with~$d \in \{1, \, 2 \}$ and the indexing of the grid points. The grid is obtained by subdividing~$(0, L)^d$ into~$N+1$ intervals along each coordinate direction. In two dimensions, only the elements in the shaded region belong to~$\Omega$. The grid points are indexed by~$k=k^0$ and~$k=(k^0, \, k^1)$ in one and two dimensions, respectively. The external boundary is subject to homogeneous Dirichlet conditions. The corresponding grid points are therefore omitted from the finite element equations. The boundary of the internal hole is subject to either homogeneous Dirichlet or Neumann conditions.    \label{fig:mu-grid}}
\end{figure}

We discretise the weak form~\eqref{eq:weak_form} with a standard finite element approach. The finite element grid is obtained by subdividing the domain~$(0, \, L)^d $ into~$N+1$ equal intervals along each coordinate direction, as shown in Figure~\ref{fig:mu-grid}.  
We assume that~$N$ is a power of~$2$. The solution~$u$ and the test function~$v$ are discretised using linear basis functions, and the diffusion coefficient~$\mu $ within each element is assumed to be constant.   The resulting finite element discretisation yields the linear system of equations 
\begin{equation}
	A \ket{u} = \ket{ f} \, , 
\end{equation}
 where~$A \in  \mathbb R^{N^d  \times N^d}$ is a sparse system matrix,~$\ket{u} \in \mathbb R^{N^d}$ is the solution array and~$\ket{f} \in \mathbb R^{N^d}$ is the force array. Although~$A$ can be singular for a multiconnected domain because of the grid points within the hole,  the corresponding singular values are suppressed during the solution stage using the quantum linear systems solver. 
 
%
\subsection{Heterogeneous system matrices \label{sec:sysMat}}
\begin{figure}[!htb]
    \centering
    \subfloat[][One-dimensional case \label{fig:1d-mu-stencil}] 
    {
    	\raisebox{-5pt}{\includegraphics[scale=0.9]{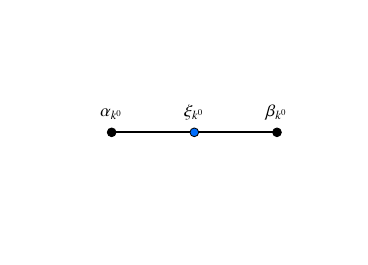} 
     }}
   \hspace{0.1\textwidth}
    \subfloat[][Two-dimensional case \label{fig:2d-mu-stencil}]
     {
     	\includegraphics[scale=0.9]{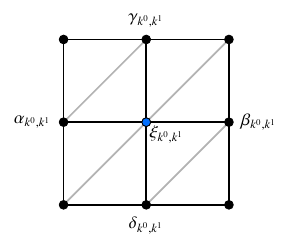} 
     }
    \caption{Local stiffness matrix stencils for a central node $k^0$ in one dimension and $(k^0, \,k^1)$ in two dimensions, respectively. The stencil weights are the entries of the system matrices corresponding to nodes $k^0$ and $(k^0, \, k^1)$.}
    \label{fig:mu-stencil}
\end{figure}

%
We focus on the simply connected domain $\Omega=(0, \, L)^d$ to expose the structure of the system matrix~$A$ of the finite element discretised problem. The boundary condition is a zero Dirichlet condition. Each row of~$A$ corresponds to a specific node in the grid, and the entries of the row depend on the diffusion coefficient in the elements sharing the node. The row entries can be conveniently represented as stencils, as shown in Figure~\ref{fig:mu-stencil}. 
%
\subsubsection{One dimension \label{sec:sysMat1D}}
%
 The three weights of the one-dimensional stencil given in Figure~\ref{fig:1d-mu-stencil} for a row corresponding to the node with the index~$k^0$ are 
\begin{equation} \label{eq:A_entries_1D}
  \xi_{k^0}  =   \frac{N+1}{L}\left(\, \mu_{k^0-1} + \mu_{k^0} \right)   \, ,  \quad \alpha_{k^0}    =  -\frac{N+1}{L}\mu_{k^0-1} \, , \quad \beta_{k^0}  =  -\frac{N+1}{L}\mu_{k^0}   \, .
\end{equation}
The assembled stencils yield the  sparse one-dimensional system matrix 
\begin{equation}
	A = 
	\begin{pmatrix*}
		\xi_0 & \beta_0 &  & &   \\
		\alpha_1 & \xi_1 & \beta_1 & &   \\
		& \alpha_2 & \ddots& \ddots &  \\
		& & \ddots & \ddots & \beta_{N-2}   \\
		& &   &  \alpha_{N-1} & \xi_{N-1}   
	\end{pmatrix*} \, . 
\end{equation}
All the omitted entries are zero. Furthermore, the system matrix of problem~$\eqref{eq:poisson}$ is usually symmetric; however, for quantum encoding, it is helpful to distinguish between the two off-diagonal bands. 

%
\subsubsection{Two dimensions \label{sec:sysMat2D}}
%
To obtain the two-dimensional stencils given in Figure~\ref{fig:2d-mu-stencil}, we subdivide each cell into two triangular finite elements. The stencil weights for the row corresponding to the node $(k^0, \, k^1)$, assembled from eight adjacent finite elements, are given by
\begin{equation} \label{eq:A_entries_2D}
\begin{aligned} 
	\xi_{k^0, k^1}  & = \mu_{k^0-1,k^1-1} + \mu_{k^0,k^1-1} + \mu_{k^0-1,k^1} + \mu_{k^0,k^1}      \, , & \\ 
	\alpha_{k^0,k^1}  & =  -\frac{1}{2}(\, \mu_{k^0-1,k^1-1} + \mu_{k^0-1,k^1}) \, , \quad  & \beta_{k^0,k^1}  =   -\frac{1}{2}(\, \mu_{k^0,k^1-1} + \mu_{k^0,k^1}) \, ,\\ 
	\gamma_{k^0,k^1}  & =  -\frac{1}{2}(\, \mu_{k^0-1,k^1} + \mu_{k^0,k^1})  \, , \quad &  \delta_{k^0, k^1}  = -\frac{1}{2}(\, \mu_{k^0-1,k^1-1} + \mu_{k^0,k^1-1})     \, .
\end{aligned}
\end{equation}
The sparse two-dimensional system matrix assembled from the stencils has the block structure
\begin{equation}
	A = 
	\begin{pmatrix*}
		C_{0,:} & D_{0,:} &   & &    \\ 
		G_{1,:} & C_{1,:} & D_{1,:} & &    \\ 
   	                     & G_{2,:} & \ddots & \ddots &  \\
	                     & & \ddots & \ddots  & D_{N-2,:}   \\
	                    &  &    &  G_{N-1,:} & C_{N-1,:}
	\end{pmatrix*} \, . 
\end{equation}
The three matrix blocks have the entries 
\begin{equation}
\begin{aligned}	
	C_{k^0,:} &= \begin{pmatrix*}
		\xi_{k^0,0} & \gamma_{k^0,0} &   & &   \\
		\delta_{k^0,1} & \xi_{k^0,1} & \gamma_{k^0,1} & &   \\
		& \delta_{k^0,2} & \ddots& \ddots &   \\
		&  & \ddots & \ddots & \gamma_{k^0,N-2}  \\
		& & &  \delta_{k^0,N-1} & \xi_{k^0,N-1}   
	\end{pmatrix*} 
	 \, , \quad 
	D_{k^0,:} = \begin{pmatrix*}[r]
		\beta_{k^0,0} & & & &     \\
		 & \beta_{k^0,1} & & &    \\
		&  &   \ddots &  &   \\
		& & &   \ddots   & \\
		& & & &  \beta_{k^0,N-1}    \\ 
	\end{pmatrix*} \, ,  \\ 
	G_{k^0,:} & = \begin{pmatrix*}
		\alpha_{k^0,0} &  & & &      \\
		 & \alpha_{k^0,1} & & &     \\
		&  & \ddots & &    \\
		& & & \ddots   & \\ 
		& & &  &  \alpha_{k^0,N-1}   
	\end{pmatrix*}  \, .
\end{aligned}
\end{equation}
%
 %
\subsection{Immersed boundary system matrices \label{sec:immersed_sys_mat}}
%
We now consider a multiconnected domain~$\Omega \subset (0, \, L)^d$ obtained by introducing an internal hole in the cuboidal domain~$(0, \, L)^d$. The corresponding system matrix is obtained by modifying the system matrix of the previously introduced heterogeneous problem. As before, we assume that the external domain boundary is homogeneous Dirichlet, and the internal hole boundary is either homogeneous Dirichlet or Neumann.  

The geometry of the hole is described by an integer-valued level-set function~$\phi (k)$, where~$k=k^0$ in one dimension and~$k=(k^0, \, k^1)$ in two dimensions. The restriction to integer-valued functions enables an efficient quantum implementation and can be relaxed, if necessary, as discussed later. The hole boundary is approximated along the lines of the Cartesian grid avoiding cut elements. The level-set function defines the three disjoint index sets 
\begin{equation}
    \mathcal{I} = \{k \in \mathbb N \forwhich \phi(k) > 0\} \, , \qquad
    \mathcal{B} = \{k \in \mathbb N \forwhich  \phi(k) = 0\} \, , \qquad
    \mathcal{H} = \{k  \in \mathbb N \forwhich \phi(k) < 0\} \, , 
\end{equation}
where~$k \in \{ 0, \, \dotsc, \,  N-1\}^d$. Here,~$\mathcal{I}$ denotes the set of all interior grid points,~$\mathcal{B}$ the set of all grid points on the hole boundary, and~$\mathcal{H}$  the set of all grid points inside the hole. 

We enforce homogeneous Dirichlet boundary conditions by setting the rows and columns of the heterogeneous stiffness matrix~$A$ associated with the grid points $k \in \mathcal{H}$ to zero. To achieve this,~$A$ is pre- and postmultiplied by a diagonal matrix~$D_\chi$. The entries of~$D_\chi$ are determined with the help of the indicator function 
\begin{equation}\label{eq:indfunc}
	\chi (k) = 
	\begin{cases} 
		1 \, ,\quad & \text{if } k \in \mathcal I \cup \mathcal B  \, , \\ 
		0 \, ,&\text{if } k \in \mathcal H  \, .
	\end{cases}
\end{equation}
 Accordingly, the immersed system matrix~$A_D$ of the Dirichlet problem is defined as 
\begin{equation}\label{eq:A_Dirichlet}
    A_D = D_\chi A D_\chi \, .
\end{equation}
For illustration, consider the one-dimensional domain with a single hole with half-width~$L/(N+1)$  centred at the grid point with index~$c$ such that~$\mathcal H = \{c\}$ and~$\mathcal B = \{ c-1, \, c+1 \}$. The respective system matrix~$A_D$  is determined as 
\begin{equation} \label{eq:A_Dirichlet_matrixForm}
\setlength{\arraycolsep}{2.0pt}
  \begin{aligned}
     \begin{pmatrix*}
          & \ddots & \ddots & \ddots &   &   &   &   \\ 
          &   & \alpha_{c-1} & \xi_{c-1} & 0  &   &   &   \\ 
          &   &   &  0 & 0 &  0 &   &   \\ 
          &   &   &   & 0  & \xi_{c+1} & \beta_{c+1} &   \\ 
          &   &   &   &   & \ddots & \ddots & \ddots
    \end{pmatrix*}
    & = 
    \begin{pmatrix*}
          \ddots &   &   &   & \\ 
          & 1 &   &   & \\ 
          &   & 0 &   & \\ 
          &   &   & 1 &  \\ 
          &   &   &   & \ddots
    \end{pmatrix*}
    \begin{pmatrix*}
          & \ddots  & \ddots & \ddots &   &   &   &   \\ 
          &   & \alpha_{c-1} & \xi_{c-1} & \beta_{c-1} &   &   &   \\ 
          &   &   & \alpha_c & \xi_c & \beta_c &   &   \\ 
          &   &   &   & \alpha_{c+1} & \xi_{c+1} & \beta_{c+1} &   \\ 
          &   &   &   &   &  \ddots & \ddots & \ddots
    \end{pmatrix*}
    \begin{pmatrix*}
          \ddots &   &   &   & \\ 
          & 1 &   &   & \\ 
          &   & 0 &   & \\ 
          &   &   & 1 &  \\ 
          &   &   &   & \ddots
    \end{pmatrix*} \, .
\end{aligned}
\end{equation}
Hence, the entries of~$D_{\chi}$ are~$\chi(k)$.

We adopt a similar approach to impose homogeneous Neumann boundary conditions on the boundary of the internal hole. The corresponding system matrix~$A_N$ is given by
\begin{equation} \label{eq:A_Neumann}
	A_N = D_\chi A D_\chi - D_{N} \left ( \frac { I \circ A }{2 d}  \right ) \, .
\end{equation}
Here, the Hadamard product~$ I \circ A $ yields the diagonal entries of the system matrix~$A$, and~$2d$ is the number of nearest neighbours of each grid point in a $d$-dimensional Cartesian grid. The diagonal matrix~$D_N$ has non-zero entries only on rows corresponding to the boundary grid points~$k \in \mathcal{B}$. The value of each entry is equal to the number of nearest neighbours of~$k$ that are in the set~$\mathcal H$. In one dimension, the number of such grid points can be only $1$, and in two dimensions, it can be $1$, $2$, or $3$ depending on the boundary geometry.

%
\section{Encoding of heterogeneous system matrices \label{sec:sysmat}}
%
Next, we consider the block-encoding of the one- and two-dimensional system matrices posed on simply connected domains with zero Dirichlet boundary conditions.  The non-unitary system matrix is, after suitable rescaling, represented as the upper-left sub-block of a larger matrix. The rescaling is necessary to ensure that all matrix entries have magnitude less than or equal to one. For reference, the block-encoding of the circulant system matrices of one-dimensional periodic problems with constant diffusion coefficients is described in~\ref{app:laplacian}.
%
%
\subsection{One dimension}
%
The action of the scaled tridiagonal matrix~$ A/ \| A\|_{\text{max}}  \in \mathbb R^{N\times N}$ on a vector~$\ket{v}  \in \mathbb R^{N}$ can be expressed as 
\begin{equation} \label{eq:AdecomposeDiag1D}
\begin{aligned} 
	 \frac{A}{\| A\|_{\text{max}}} \ket{v} 
	& = 
	\begin{pmatrix*}
		\tilde{\xi}_0 &  &  &  &   \\
		& \tilde{\xi}_1 & & &   \\
		&  & \ddots & &  \\
		& &  & \ddots &    \\
		& &  &  & \tilde{\xi}_{N-1}   
	\end{pmatrix*} 
	\begin{pmatrix*}
		v_0  \\
		v_1  \\
		 \vdots    \\
		 v_{N-2} \\
		v_{N-1}    
	\end{pmatrix*}    + 
	\begin{pmatrix*}
		\tilde{\alpha}_0 &  &  &  &   \\
		& \tilde{\alpha}_1 & & &   \\
		&  & \ddots & &  \\
		& &  & \ddots &    \\
		& &  &  & \tilde{\alpha}_{N-1}   
	\end{pmatrix*} 
		\begin{pmatrix*}
		0 \\
		v_0  \\
		\vdots    \\
		v_{N-3} \\
		v_{N-2}    
	\end{pmatrix*} \\ & +
	\begin{pmatrix*}
		\tilde{\beta}_0 &  &  &  &   \\
		& \tilde{\beta}_1 & & &   \\
		&  & \ddots & &  \\
		& &  & \ddots &    \\
		& &  &  & \tilde{\beta}_{N-1}   
	\end{pmatrix*} 
		\begin{pmatrix*}
		v_1  \\
		v_2  \\
		 \vdots    \\ 
		 v_{N-1}\\
		0    
	\end{pmatrix*}     \, , 
\end{aligned}
\end{equation}
with the maximum norm~$\| A\|_{\text{max}}$ and the rescaled entries 
\begin{equation} \label{eq:A_entries_1D_rescaled}
	\tilde{\xi}_{k^0} = \frac{\xi_{k^0}}{ \| A\|_{\text{max}}}  \, ,  \quad \tilde{\alpha}_{k^0} = \frac{ \alpha_{k^0} }{ \| A\|_{\text{max}}} \, , \quad  \tilde{\beta}_{k^0} = \frac{ \beta_{k^0} }{ \| A\|_{\text{max}}} \, . 
\end{equation}
Hence, $ (A/ \| A\|_{\text{max}} ) \ket{v}$ is the sum of three diagonal matrices applied to the shifted versions of $\ket{v}$  (with zero padding).  We define the two shift operators 
\begin{equation} \label{eq:shiftOperators}
	 S_{+} 
		\begin{pmatrix*}
		v_0  \\
		v_1  \\
		 \vdots    \\
		 v_{N-2} \\
		v_{N-1}    
	\end{pmatrix*} 
	= \begin{pmatrix*} 0 \\ v_0  \\ \vdots    \\v_{N-3} \\v_{N-2}    \end{pmatrix*}
	\, , \qquad    S_{-}   
		\begin{pmatrix*}
		v_0  \\
		v_1  \\
		 \vdots    \\
		 v_{N-2} \\
		v_{N-1}    
	\end{pmatrix*}   = 
	\begin{pmatrix*} v_1  \\ v_2  \\  \vdots    \\  v_{N-1}\\ 0     \end{pmatrix*}  \, . 
\end{equation}

\begin{figure*}[!tbp]
    \centering
    \scalebox{1}{
    	\includegraphics[scale=1]{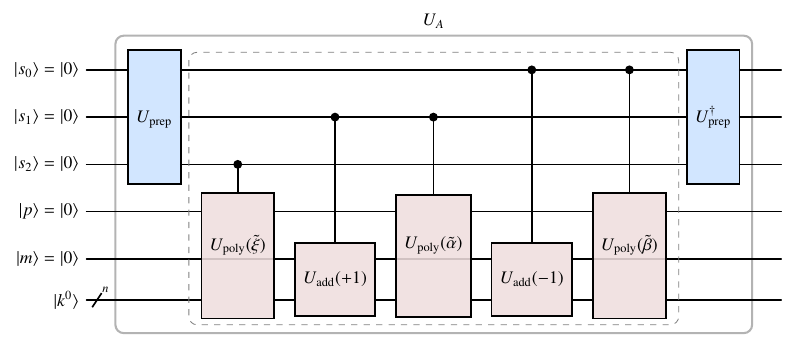} 
    }
    \caption{Block encoding~$U_A$ of the scaled system matrix~$A/ \| A\|_{\text{max}}$ of the one-dimensional problem.  The discretisation consists of~$N=2^n$ internal grid points.  The three $U_{\text{poly}}$ gates encode the main, sub- and superdiagonal of~$A/ \| A\|_{\text{max}}$. The two~$U_{\text{add}}$ gates implement the up and down shifts of the components of the input~$\ket {k^0}$. The sum of the main, sub- and superdiagonal entries is computed using the LCU construction with the aid of~$U_{\text{prep}}$ and its inverse~$U_{\text{prep}}^\dagger$. In the output, the relevant components correspond to basis states~$\ket{0}^{\otimes 3}\ket{0}\ket{0}\ket{k^0}$. \label{fig:1DhetroCirc}}
\end{figure*}

\begin{figure}[!tbp]
\centering
	\subfloat[][ ]{ \label{fig:poly_and_add_a}
	\scalebox{1}{
    	\includegraphics[scale=1]{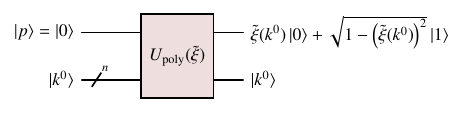}  }	
	}
	\hspace{0.04\textwidth}
	\subfloat[][ ]{ \label{fig:poly_and_add_b}
	\scalebox{1}{
    	\includegraphics[scale=1]{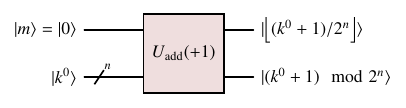} }	
	}
\caption{Polynomial encoding and modular addition gates used in block encoding.  In~(a) the amplitude of the output basis state $\ket 0 \ket{k^0}$ is the sought function value~$\tilde{\xi}(k^0)$.  In~(b), the top qubit acts as a carry (overflow) bit in modular arithmetic. It ensures that the output amplitudes are padded with zeros.
\label{fig:poly_and_add}}
\end{figure}

The circuit implementation of the decomposition~\eqref{eq:AdecomposeDiag1D} is depicted in Figure~\ref{fig:1DhetroCirc}. 
The circuit encodes the scaled matrix~$A/ \| A\|_{\text{max}}$, with entries bounded in magnitude by~$1$. As defined in~\eqref{eq:A_entries_1D_rescaled} and~\eqref{eq:A_entries_1D}, the diagonal entries~$\tilde{\xi}_{k^0}$, $\tilde{\alpha}_{k^0}$, and~$\tilde{\beta}_{k^0}$ depend on the diffusion coefficient~$\mu$ evaluated over the elements sharing the node~$k^0$. If $\mu$ is a polynomial, then the corresponding functions $\tilde{\xi}(k^0)$, $\tilde{\alpha}(k^0)$, and~$\tilde{\beta}(k^0)$ are polynomials of the same degree; otherwise, they are approximated classically by polynomials. These polynomials are then quantum encoded using standard techniques. For instance, the encoding of~$\tilde{\xi}(k^0)$ using the unitary~$U_{\text{poly}} (\tilde{\xi})$ yields
\begin{equation} \label{eq:Upoly}
	U_{\text{poly}}( \tilde{\xi} ) \colon \ket 0  \ket{k^0} \mapsto \tilde{\xi}(k^0) \ket 0 \ket{k^0} +  \sqrt{1-  \left (\tilde{\xi}(k^0) \right)^2}  \ket 1 \ket{k^0} \, . 
\end{equation}
 The input consists of the single ancilla qubit in state~$\ket 0$ and the $n$-qubit register~$\ket {k^0}$, where~$n = \log_2 N $.  The register~$\ket {k^0}$ provides a labelling of the components of the vector to which the scaled matrix~$A/ \| A\|_{\text{max}}$ is applied. The unitaries~$U_{\text {add}}(+1)$ and~$U_{\text {add}}(-1)$ implement the non-unitary shift operators~$S_{+}$ and~$S_{-}$, respectively. For instance, the unitary $U_{\text {add}}(+1)$ is defined as
\begin{equation} \label{eq:Uadd}
	U_{\text {add}}(+1) : \ket 0 \ket {k^0}  \mapsto  \ket{\left\lfloor (k^0+1)/N \right\rfloor} \ket { (k^0+1) \mod N}  \, , 
\end{equation}
where the first qubit serves as a carry bit, with the floor operator representing overflow. $U_\text{add}$ is realised as modular (cyclic) addition, and the first qubit is used to pad the shifted vector with zero. This is achieved by mapping the zero-amplitude input component~$\ket{1}\ket{N-1}$  to~$\ket {0} \ket{0}^{\otimes n}$.  The unitary~$U_{\text {add}}(-1)$ is defined analogously. See Figure~\ref{fig:poly_and_add} for the circuit representations of~\eqref{eq:Upoly} and~\eqref{eq:Uadd}. Further details on the implementation  of~$U_{\text{poly}}$  and~$U_{\text {add}}$ can be found in~\cite{liu2024towards} and~\ref{app:qc-addition}, respectively. 

The sum of the three components in decomposition~\eqref{eq:AdecomposeDiag1D} is computed using the LCU construction. To this end, the unitary~$U_{\text{prep}}$ generates the unary, or the~$W$, state
\begin{equation}
	U_{\text{prep}} \colon \ket 0 \ket 0 \ket 0 \mapsto \frac{1}{\sqrt 3} \left ( \ket 0 \ket 0 \ket 1 + \ket 0 \ket 1\ket 0 + \ket 1 \ket 0 \ket 0 \right ) \, ,
\end{equation}
where exactly one of the qubits is in state~$\ket 1$. The controlled~$U_{\text{poly}}$ and~$U_{\text{add}}$ gates in the circuit constitute the selection layer~$U_{\text{select}}$ in the LCU construction. The resulting circuit, referred to as the~$U_A$ gate, is a block encoding of the scaled stiffness matrix~$A/ \| A\|_{\text{max}}$ and implements the mapping
\begin{equation}
	U_A \colon \sum_{k^0=0}^{N-1} v_{k^0} {\ket 0}^{\otimes 3} \ket{0}\ket{0} \ket{k^0}  \mapsto \sum_{k^0=0}^{N-1} \frac{1}{3 \| A\|_{\text{max}} } (A v)_{k^0} {\ket 0}^{\otimes 3} \ket{0}\ket{0} \ket{k^0} + \dotsc  \, , 
\end{equation}
which can be interpreted as 
\begin{equation} \label{eq:Av_blockMatrixForm}
	\begin{pmatrix} \dfrac{A}{\sigma}& *  \\ * & * \end{pmatrix} \begin{pmatrix} \ket v \\ 0  \end{pmatrix} = \begin{pmatrix} \ \dfrac{\ket {Av}}{\sigma}   \\ *  \end{pmatrix} \, ,
\end{equation}
where~$\sigma=3 \| A\|_{\text{max}}$ is the subnormalisation factor of the introduced block-encoding. The stars indicate submatrices and vectors that are introduced by the LCU construction and ensure unitarity. The dimension of the entire vector space in~\eqref{eq:Av_blockMatrixForm} is~$2^5 N$. The relevant components correspond to the first~$N$ basis states. 

We note that other block encodings of the system matrix are feasible. In general, block encodings with a small subnormalisation factor~$\sigma$ are preferred, as the postselection probability of the first~$N$ states, i.e. where the ancillas are in state~${\ket 0}^{\otimes 3} \ket 0 \ket 0$,  is proportional to~$ \| A \ket v\|^2 / \sigma^2 $. Furthermore, the choice of the~$W$ state in the proposed LCU construction ensures that each~$U_{\text{select}}$ unitary requires only a single control. In contrast, using a uniform superposition state prepared with Hadamard gates requires multi-controlled unitaries in $U_{\text{select}}$. Consequently, the choice of the~$W$ state leads to a significant reduction in the number of~$\mathit{CNOT}$ gates.

%
%
\subsection{Two dimensions}
%
We consider the action of the scaled block tridiagonal matrix~$A /  \| A\|_{\text{max}} \in \mathbb R^{N^2\times N^2}$ on a vector~$\ket {v} \in \mathbb R^{N^2}$. The components of~$\ket{v}$ are indexed by the multi-index~$(k^0, \, k^1)$, with~$k^0 \in \{ 0, \, 1, \, \dotsc, \, N-1 \}$ and~$k^1 \in \{ 0, \, 1, \, \dotsc, \, N-1 \}$. The matrix-vector product can be written as 
\begin{equation} \label{eq:AdecomposeDiag2D}
\begin{aligned}
	\frac{A}{\| A\|_{\text{max}}} \ket{v} &=
	\begin{pmatrix*}
		\tilde C_{0,:} &  &   & &    \\ 
			    & \tilde C_{1,:} &  & &    \\ 
   	                     & & \ddots &  &  \\
	                     & &  & \ddots  &   \\
	                    &  &    &   & \tilde C_{N-1,:}
	\end{pmatrix*} 
	\begin{pmatrix}
		v_{0, :} \\  v_{1, :} \\  \vdots \\  v_{N-2, :} \\ v_{N-1, :}  
	\end{pmatrix}   
	+	\begin{pmatrix*}
		\tilde D_{0,:} &  &   & &    \\ 
			    & \tilde D_{1,:} &  & &    \\ 
   	                     & & \ddots &  &  \\
	                     & &  & \ddots  &   \\
	                    &  &    &   & \tilde D_{N-1,:}
	\end{pmatrix*} 
	\begin{pmatrix}
		v_{1, :} \\  v_{2, :} \\  \vdots \\  v_{N-1, :} \\ 0  
	\end{pmatrix}  \\
	& +	\begin{pmatrix*}
		\tilde G_{0,:} &  &   & &    \\ 
			    & \tilde G_{1,:} &  & &    \\ 
   	                     & & \ddots &  &  \\
	                     & &  & \ddots  &   \\
	                    &  &    &   & \tilde G_{N-1,:}
	\end{pmatrix*} 
	\begin{pmatrix}
		0 \\  v_{0, :} \\  \vdots \\  v_{N-3, :} \\    v_{N-2, :} 
	\end{pmatrix}  \, , 
\end{aligned}
\end{equation}
with the rescaled block matrices 
\begin{equation}
	\tilde C = \frac{C}{\| A\|_{\text{max}}} \,  , \quad \tilde D = \frac{D}{\| A\|_{\text{max}}} \,  , \quad \tilde G = \frac{G}{\| A\|_{\text{max}}} \, .
\end{equation}
 The colon in indices in~\eqref{eq:AdecomposeDiag2D} indicates all entries along the corresponding dimension. The action of each block~$C_{k^0,:}$ is implemented analogously to the one-dimensional case using the two-dimensional shift operators~$I \otimes S_+$ and~$I \otimes S_-$, where $I \in \mathbb R^{N\times N }$ denotes the identity matrix, and $S_{+}$ and $S_{-}$ are defined in~\eqref{eq:shiftOperators}.  The shift operators for the second and third terms   in~\eqref{eq:AdecomposeDiag2D}  are~$S_{-} \otimes I$ and~$S_{+} \otimes I$, respectively,  and act on the first index of the multi-index.  

\begin{figure*}[!tbp]
    \centering
    \includegraphics[scale=0.93]{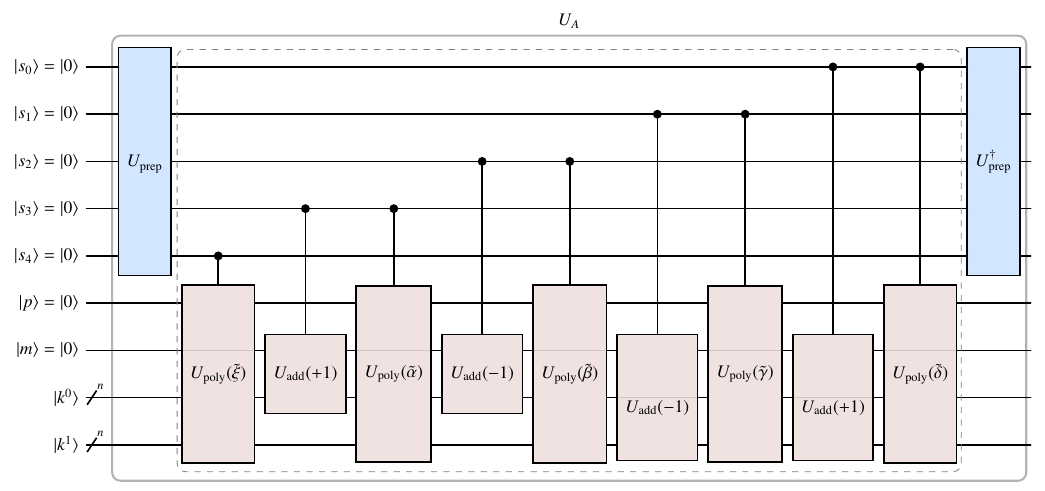} 
    \caption{Block encoding~$U_A$ of the system matrix~$A/\| A\|_{\text{max}}$ of the two-dimensional problem. The discretisation consists of \mbox{$(N+1)^2$} cells and \mbox{$N^2 =2^{2n}$} internal grid points.  The five $U_{\text{poly}}$ gates encode the bivariate polynomials~$\tilde \xi(k^0,k^1)$, $\tilde \alpha(k^0,k^1)$, $\tilde \beta(k^0,k^1)$, $\tilde \gamma(k^0,k^1)$ and $\tilde \delta(k^0,k^1)$. The four~$U_{\text{add}}$ gates implement the up and down shifts of the indices in the multi-index~$(k^0, k^1)$. The sum of these contributions is computed using the LCU construction with the aid of~$U_{\text{prep}}$ and its inverse~$U_{\text{prep}}^\dagger$. In the output, the relevant components correspond to basis states~$\ket {0}^{\otimes 5} \ket{0}\ket{0}\ket{k^0}\ket{k^1}$. \label{fig:2DhetroCirc}}
\end{figure*}

The circuit implementation of the decomposition~\eqref{eq:AdecomposeDiag2D} is shown in Figure~\ref{fig:2DhetroCirc}.  
There are in total five components in the decomposition~\eqref{eq:AdecomposeDiag2D}, and their sum is evaluated using LCU. Specifically, the gate~$U_{\text{prep}}$ creates the five-qubit $W$ state
\begin{equation}
	U_{\text{prep}} \colon \ket 0^{\otimes 5} \mapsto \frac{1}{\sqrt 5} \left ( \ket {00001} + \ket {00010} + \ket {00100} + \ket {01000} + \ket {10000}   \right ) \, .
\end{equation}
The select layer contains the five controlled~$U_{\text{poly}}$ gates encoding the bivariate functions~$\tilde{\xi}(k^0, k^1) $, $\tilde{\alpha}(k^0, k^1)$, $\tilde{\beta}(k^0, k^1)$, $\tilde{\gamma}(k^0, k^1)$ and $\tilde{\delta}(k^0, k^1)$, which are scaled versions of the corresponding functions defined in~\eqref{eq:A_entries_2D}.  Each function depends on the diffusion coefficient~$\mu$ evaluated on the elements sharing the node~$(k^0, \,  k^1)$.  If the diffusion coefficient~$\mu$ is a polynomial, each of the bivariate functions is a polynomial; otherwise, they are first approximated classically by bivariate polynomials and then quantum encoded. For instance, the encoding of~$\tilde{\xi}(k^0, k^1)$ using the unitary~$U_{\text{poly}}(\tilde \xi)$ reads
\begin{equation} \label{eq:Upoly_2d}
	U_{\text{poly}}( \tilde{\xi} ) \colon \ket 0  \ket{k^0} \ket{k^1} \mapsto \tilde{\xi}(k^0, k^1) \ket 0 \ket{k^0}  \ket{k^1}  +  \sqrt{1-  \left(\tilde{\xi}(k^0, k^1) \right)^2}  \ket 1 \ket{k^0} \ket{k^1} \, . 
\end{equation}
The implementation of the four controlled shift operators~$I \otimes S_{+}$, $I \otimes S_{-}$, $S_{+} \otimes I$, $S_{-} \otimes I$ follows the one-dimensional implementation. 

The entire circuit implements the mapping 
\begin{equation}
	\sum_{k^0, \, k^1=0}^{N-1} v_{k^0 k^1} {\ket 0}^{\otimes 5} \ket{0}\ket{0} \ket{k^0} \ket{k^1} \mapsto \sum_{k^0, \, k^1=0}^{N-1} \frac{1}{5 \| A\|_{\text{max}} } (A v)_{k^0 k^1} {\ket 0}^{\otimes 5} \ket{0}\ket{0} \ket{k^0} \ket{k^1} + \dotsc  \, , 
\end{equation}
and has the subnormalisation factor~$\sigma = 5 \| A \|_{\text{max}}$, cf.~\eqref{eq:Av_blockMatrixForm}.

 %
 \section{Encoding of immersed system matrices \label{sec:immersed}}
 %
In this section, we consider the block-encoding of system matrices for problems posed in multiconnected domains with an internal hole. As before, the external domain boundary is homogeneous Dirichlet, and the internal hole boundary is either homogeneous Dirichlet or Neumann. The hole geometry is defined by a level-set function which is in turn used to define a binary indicator function.  The block-encoded matrix is obtained by modifying the block-encoded system matrix of the heterogeneous problem using the indicator function. 
%
\subsection{Quantum level-set and indicator functions \label{sec:level-setAndIndicator}}
%
Any level-set or indicator function can, in principle, be first classically approximated using a polynomial and then quantum encoded. However, this would require extremely high polynomial degrees and is prone to oscillatory Gibbs phenomena. Therefore, we represent both functions using basis encoding and compute them directly using quantum arithmetic. 

In basis encoding, integers are represented in binary using computational basis states. Briefly, e.g., the integer $k^0=6$ (in binary $110$) is represented by the basis vector~$\ket {1 1 0}$. This makes it possible to compute, e.g., the sum~$k^0+3$ using quantum arithmetic: 
\begin{equation}
	 U_{\text{add}} \colon \overset{m}{\ket{0}} \underbrace{\ket {1 1 0}}_{k^0=6} \underbrace{\ket {0 1 1}}_{=3} \mapsto \underbrace{ \overset{m}{\ket 1} \ket {0  0 1}}_{k^0+3=9} \underbrace{\ket {0   1 1}}_{=3}  \, .
\end{equation}
Here, the most significant qubit $\ket m$ of the output register serves as a carry (overflow) bit for addition, and a borrow (underflow) bit for subtraction. The least significant three qubits~$\ket{011}$ are maintained to ensure reversibility, as all quantum operations must be unitary. A possible implementation of~$U{_\text{add}}$ using a QFT-based modular arithmetic approach is given in~\ref{app:qq-addition}. In this paper, we consider only integer-valued level-set functions. This restriction can be relaxed by allowing fixed-point numbers rather than integers in quantum arithmetic as discussed in \ref{app:floating_point}.
 
 We next introduce our basic approach to computing the level-set functions, using illustrative examples in one and two dimensions.  

%
\subsubsection{One dimension\label{sec:oneDlsetAndIndic}}
\begin{figure}[!tbp]
    \centering
    	\subfloat[][Domain with a hole and grid points \label{fig:1d-domain_immersed}] {
		\includegraphics[scale=0.95]{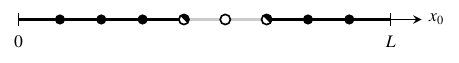} 
	}
   	\hfill
	\subfloat[][Level-set and indicator functions  \label{fig:1d-signed_distance}] {
		\includegraphics[scale=0.95]{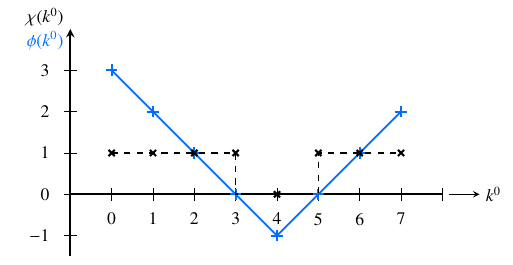} 	
        }
    \caption{One-dimensional domain with an internal hole discretised with~$N+1=9$ elements. The hole geometry is described by a level-set function~$\phi(k^0)$, which is used to compute the indicator function~$\chi(k^0)$. The grid points corresponding to the external Dirichlet boundaries are omitted from the finite element equations and the level-set and indicator function computations. \label{fig:domain_w_hole}}
\end{figure}

%
As an example, we consider the one-dimensional domain discretised with~$N=8$ internal grid points as shown in Figure~\ref{fig:domain_w_hole}. 
The grid points at the outer domain boundaries~$x_0=0$ and~$x_0=L$ are omitted because of the homogeneous Dirichlet boundary conditions at the outer boundaries.  The geometry of the hole, centred at the grid point~$c=4$ and with half-width~$b=1$, is composed of the level-set functions of two half-lines 
\begin{equation}\label{eq:sgndist_1d_1}
	\phi (k^0) = \max \{ -k^0 +c ,  \,  k^0 -c  \} - b \, .
\end{equation}
Expressing the maximum in terms of the absolute value yields   
\begin{equation}\label{eq:sgndist_1d_2}
	\phi (k^0) = | k^0 - c | - b \, .
\end{equation}
The index sets for the interior, hole boundary and hole grid points are \mbox{$\mathcal{I}= \{0, \, 1, \, 2, \, 6, \, 7 \}$}, \mbox{$\mathcal{B} =\{3, \, 5\}$} and \mbox{$\mathcal{H}=\{4\}$}, respectively. 

\begin{figure}[!tbp]
   \centering
	\subfloat[][Level-set unitary~$U_{\phi}$]{ \label{fig:1d_lvlset_circ}
		\includegraphics[scale=1]{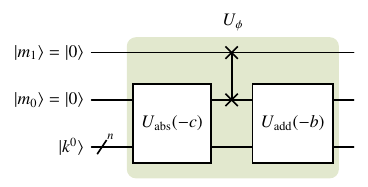} 	
	}
	\hspace{0.1\textwidth}
	\subfloat[][Indicator unitary~$U_\chi$]{ \label{fig:1d_indicator_circ}
		\includegraphics[scale=1]{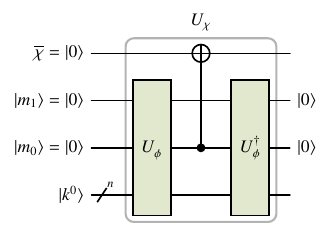} 
	}
	\caption{Computing the level-set and indicator functions~$\phi(k^0)$ and~$\chi (k^0)$ for a one-dimensional domain. (a) The level-set unitary~$U_\phi$ encodes a hole with  half-width~$b$ and centre~$c$.  (b) The indicator unitary~$U_\chi$ sets the state of the qubit~$\ket{\overline \chi}$ to~$\ket 1$ when~$\ket {m_0}$ is in state~$\ket 1$ for~$\phi(k^0)<0$ for all \mbox{$k^0\in \mathcal{H}$}.  \label{fig:levelset_and_indicator}}
\end{figure}

 The quantum arithmetic circuit for evaluating the level-set function~$\phi (k^0) $ is depicted in Figure~\ref{fig:1d_lvlset_circ}. The level-set function is calculated by first applying the unitary~$U_{\text{abs}}$ to a given grid index~$k^0$ to obtain the absolute value~$|k^0 - c|$, where the qubit~$\ket{m_0}$ is an ancilla that serves as a borrow bit. Before applying~$U_{\text{add}}$, it is swapped with the initialised ancilla~$\ket {m_1}$. The unitary~$U_{\text{add}}$ subtracts the half-width~$b$ from the absolute value~$ | k^0 -c |$. The ancilla~$\ket {m_0} $ is the borrow bit and is in state~$\ket 1$ when~$\phi(k^0) <  0$ and in state~$\ket 0$ when~$\phi(k^0) \ge 0$. 

As shown in the circuit in Figure~\ref{fig:1d_indicator_circ},  the state of the ancilla~$\ket {m_0}$ is used to set the state of the indicator qubit~${\ket{\overline \chi}}$. The indicator qubit~$\overline \chi$ is in state~$\ket 1$ when~$k^0 \in \mathcal H$ and in state~$\ket 0$ otherwise. Hence, the indicator qubit~$\ket{\overline \chi}$ encodes the function~$\overline \chi(k^0) = 1 -\chi(k^0)$. 
 After $\ket {\overline \chi}$ is set, the two ancillary qubits~$\ket {m_0}$ and~$\ket {m_1}$ are reset to their original states~$\ket 0$ by uncomputing, i.e. by applying~$U^\dagger_\phi$.  This is essential because the level-set and indicator functions must be recomputed every time the block-encoded system matrix is invoked, e.g., when solving the linear systems of equations. 
  
 The constructions of the quantum arithmetic unitaries~$U_{\text{add}}$ and~$U_{\text{abs}}$ are discussed in~\ref{app:qc-addition} and~\ref{app:absolute}, respectively. Both can be implemented using QFT-based arithmetic operations, although other implementations are available in quantum programming libraries, such as Qiskit~\cite{qiskit2024} and PennyLane~\cite{bergholm2018pennylane}. Furthermore, it bears emphasis that these quantum circuits can simultaneously compute the level-set and indicator functions for all grid points when a uniform superposition state is provided to the register $\ket{k^0}$.
 %
 \subsubsection{Two dimensions\label{sec:twoDlsetAndIndic}}
\begin{figure}[!tbp]
    \centering
    	\subfloat[][Level-set circuit~$U_\phi$ \label{fig:2d_signed_distance_circ}] {
		\includegraphics[scale=1]{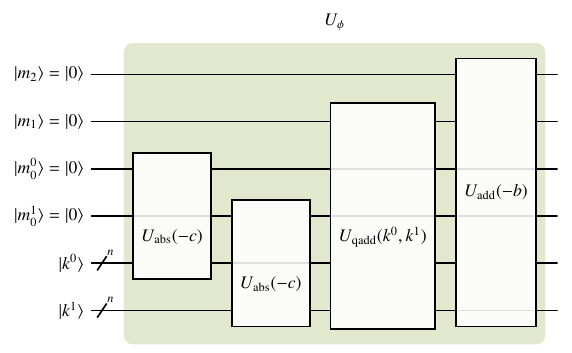} 
	}
   	\hfill
	\subfloat[][Indicator function~$\chi(k^0,k^1)$  \label{fig:2d-indicator_func}] {
		\includegraphics[scale=1]{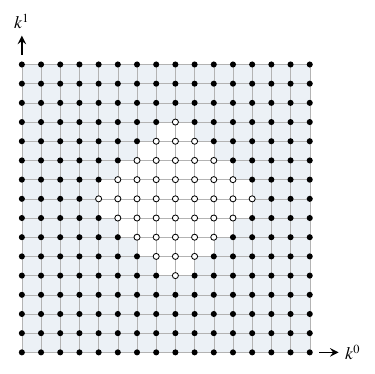} 	
        }
    	\caption{Computing the level-set function~$\phi(k^0, k^1)$ of a two-dimensional problem with a diamond-shaped hole with a half-width~$b$ centred at grid point~$(c, \, c)$. (a) Quantum circuit for evaluating~$\phi(k^0, k^1)$ by computing the signed Manhattan distance to the hole centre. (b) The corresponding indicator function~$\chi(k^0, k^1)$ evaluated at the grid points. At full dots, the indicator function has the value $1$, and at the empty dots, the value~$0$. The full dots indicate the grid points in the index set~$\mathcal I \cup \mathcal B$ and the empty dots the grid points in the index set~$\mathcal H$. \label{fig:levelset_and_indicator_2d}}
\end{figure}

 %
As a two-dimensional example, we consider the domain with the diamond-shaped hole shown in Figure~\ref{fig:2d-indicator_func}.  Its level-set function can be composed from the level-set functions of four half-planes using their maximum, 
\begin{equation}
	\phi(k) = \max \{ (k^0-c)+(k^1-c), \, (k^0-c)-(k^1-c) , \, -(k^0-c)+(k^1-c) , \, -(k^0-c)-(k^1-c) \} - b \, .
\end{equation}
Expressing the maximum in terms of absolute values yields  
\begin{equation}\label{eq:sgndist_2d}
	\phi (k) = | k^0 - c | + | k^1 - c | - b \, .
\end{equation}
The hole is centred at the grid point with the index~$(c, \,c)$ and has the half-width~$b$. The level-set function represents the signed Manhattan distance of the grid point~$k=(k^0, \,  k^1) $ to the grid points in the boundary index set~$\mathcal{B}$.

The quantum arithmetic circuit for computing the level-set function is presented in~Figure~\ref{fig:2d_signed_distance_circ}. Similar to the one-dimensional case, the level-set function~$\phi(k)$  is computed by first applying the unitary~$U_{\text{abs}}$ to each index of the grid point~$k$. The two ancillas~$\ket{m_0^0}$ and~$\ket{m_0^1}$ are the borrow bits for the two~$U_{\text{abs}}$ unitaries. Subsequently, the absolute values~$| k^0 - c |$ and~$| k^1- c |$ are added by applying the quantum addition unitary~$U_{\text{qadd}}$, which has the carry bit~$\ket{m_1}$.  The last unitary~$U_{\text{add}}$ subtracts the half-width~$b$.   Finally, the borrow bit~$\ket{m_2}$ can be used to set the indicator bit~$\overline \chi$ analogous to the one-dimensional circuit in Figure~\ref{fig:1d_indicator_circ}. 

 The constructions of the unitaries~$U_{\text{add}}$, $U_{\text{abs}}$ and~$U_{\text{qadd}}$ are discussed in~\ref{app:qc-addition},~\ref{app:absolute} and~\ref{app:qq-addition}, respectively.

 %
 \subsection{Dirichlet problems \label{sec:immersed_sysmat_D}}
 %
As discussed in Section~\ref{sec:immersed_sys_mat}, we consider homogeneous Dirichlet conditions at the hole boundary by pre- and post-multiplying the system matrix~$A$ of the simply connected domain by the diagonal indicator matrix~$D_\chi$.  As given in~\eqref{eq:A_Dirichlet}, the system matrix of the Dirichlet problem is~\mbox{$  A_D = D_\chi A D_\chi $}.  The entries of~$D_\chi$ corresponding to the grid points in the hole, i.e. for $k \in \mathcal{H}$, are~$0$, and~$1$ otherwise. 

\begin{figure}[!tbp]
    \centering
        \centering
            \includegraphics[scale=1]{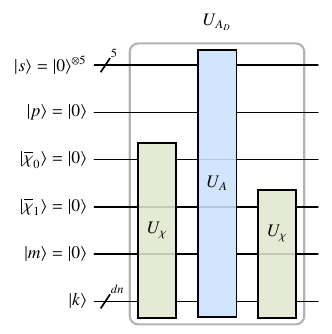}        
        \caption{Block encoding~$U _ {A_D} $ of the scaled immersed system matrix~$A_D/ \| A_D \|_{\text{max}}$ of a problem with homogeneous Dirichlet conditions at the hole boundary, where~$ A_D = D_\chi A D_\chi $.  The two $U_\chi$ unitaries block-encode the indicator matrix~$D_\chi$. The $U_A$ unitary block-encodes the system matrix of the simply connected problem.	 \label{fig:twoD_immersed_dirichlet}}
\end{figure}

 In the quantum setting, the introduced unitary~$U_{\chi}$ for the indicator function provides a block-encoding of the matrix~$D_\chi$. The unitaries~$U_A$ introduced in Section~\ref{sec:sysmat} provide block-encodings of one- and two-dimensional system matrices~$A$ for simply connected domains. The block-encoding of the system matrix~$A_D$ can be constructed from the block-encodings~$U_{\chi}$ and~$U_A$.  Products of block-encoded matrices are straightforward to compute when separate ancilla registers are used for their block-encodings~\cite{gilyen2019quantum}; see also~\ref{app:product_encoding}. The quantum circuit for the unitary~$U_{A_D}$, which provides a block-encoding for~$A_D$, is depicted in Figure~\ref{fig:twoD_immersed_dirichlet}. At the circuit output, the relevant components correspond to basis states for which all registers other than the~$\ket k$ register are in the state~$\ket 0$.
%
  
 %
 \subsection{Neumann problems \label{sec:immersed_sysmat_N}}
 %
 %
We obtain the immersed system matrix~$A_N$  for the problem with homogeneous Neumann boundary conditions on the hole boundary by modifying the Dirichlet system matrix~$A_D$. Specifically, the diagonal entry associated with each boundary grid point is reduced in proportion to the number of its nearest-neighbour grid points inside the hole.  According to~\eqref{eq:A_Neumann}, the system matrix of the Neumann problem is~\mbox{$  A_N = A_D - D_{N} (  I \circ A )/ (2 d) $}.  The diagonal matrix~$D_N$ has non-zero entries only at grid points~$k$ on the hole boundary, i.e.~$k \in \mathcal{B}$. For each such grid point, the corresponding entry is equal to the number of its nearest neighbours belonging to~$\mathcal{H}$.

\begin{figure}[!tbp]
    \centering
    \includegraphics[scale=1]
        {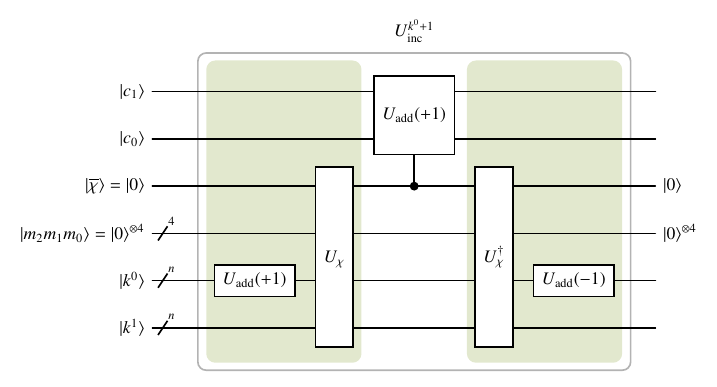}
	\caption{Quantum circuit for determining whether the right neighbour~$(k^0+1, \, k^1)$ of the grid point $(k^0, \, k^1)$ lies inside the hole. When \mbox{$(k^0+1, \, k^1) \in \mathcal{H}$}, the indicator unitary $U_\chi$ sets $\ket{\overline \chi}$ to $\ket 1$, and the counter register~$\ket{c_1c_0}$ is incremented by one. \label{fig:2DImmNCounterCirc}}\end{figure}

The entries of the diagonal matrix~$D_N$ are therefore determined by counting for each boundary grid point its nearest neighbours inside the hole. This can be accomplished by inspecting the two nearest neighbours along each coordinate direction and incrementing a counter whenever a neighbour belongs to~$\mathcal H$. As an example, Figure~\ref{fig:2DImmNCounterCirc} depicts  the circuit~$U_{\text{inc}}^{k^0+1}$, which determines whether the nearest neighbour~$(k^0+1, \, k^1)$  of the grid point~$(k^0, k^1)$ lies inside the hole.  The circuit consists of three substeps, with the first and third substeps marked by shaded regions. In the first substep, the grid-point index is first incremented by applying the unitary~$U_{\text {add}}$. Subsequently, the indicator unitary~$U_{\chi}$ introduced in Section~\ref{sec:level-setAndIndicator} is applied to determine whether the neighbouring grid point lies inside the hole. In the second substep, the counter register~$\ket {c_1c_0}$ is incremented conditional on~$\ket{\overline \chi} = \ket 1$. Finally, in the third substep, the indicator and ancilla registers are uncomputed by applying~$U_{\chi}^\dagger$, and the original grid point index is restored by applying~$U_{\text{add}}(-1)$. 

\begin{figure}[!tbp]
    \centering
    \includegraphics[scale=0.90, trim=10pt 5pt 10pt 5pt,clip]{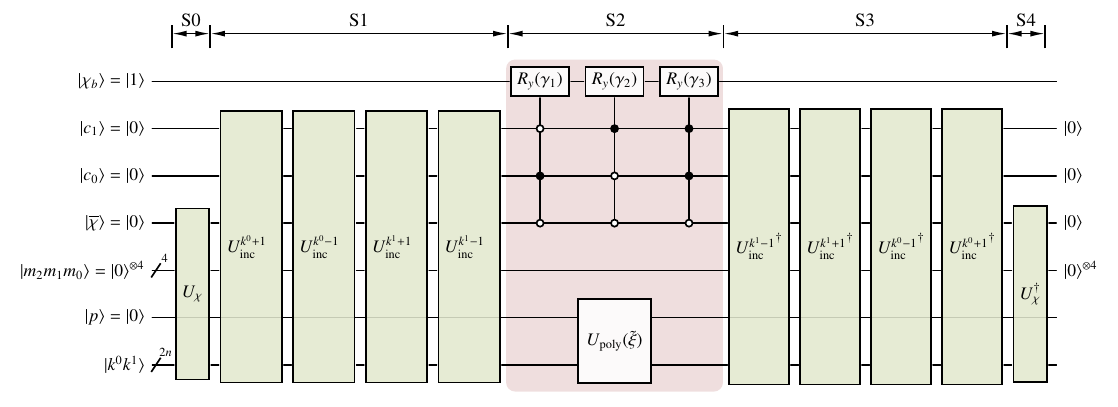}
    \caption{Block encoding of the diagonal matrix $D_{N} (  I \circ A/\| A \|_{\text{max}} )/ 4$. The entries of~$D_N$ are determined by counting the nearest neighbours of the grid point~$(k^0, \, k^1)$ that belong to~$\mathcal{H}$. The entries~$\tilde \xi$ of~$ I \circ A/\| A \|_{\text{max}} $ are scaled by applying the~$R_Y$ gates with rotation angles that depend on the number of nearest neighbours in~$\mathcal{H}$.}
    \label{fig:2DImmNCountBoundaryCirc}
\end{figure}

The diagonal correction matrix~$D_{N} (  I \circ A/\| A \|_{\text{max}} )/ (2 d)$ is block-encoded using the circuit shown in Figure~\ref{fig:2DImmNCountBoundaryCirc}, which consists of the five substeps denoted by S0, S1, S2, S3 and S4. In S0,  the indicator unitary~$U_{\chi}$ is applied to examine the location of grid point~$(k^0, \, k^1)$ in relation to the hole. The state of qubit~$\ket{\overline \chi}$ is set to~$\ket 1$ when $(k^0, \, k^1) \in \mathcal {H}$, otherwise it remains~$\ket 0$.  In S1, the number of the neighbouring grid points belonging to~$\mathcal H$ is counted by applying the unitaries~$U_{\text{inc}}^{k^0+1}$, $U_{\text{inc}}^{k^0-1}$,  $U_{\text{inc}}^{k^1+1}$ and $U_{\text{inc}}^{k^1-1}$ and incrementing the counter register~$\ket {c_1} \ket {c_0}$. Subsequently, in S2 the diagonal entry~$ I \circ A$ is evaluated by applying the unitary~$U_{\text{poly}} (\tilde \xi)$ defined in~\eqref{eq:Upoly_2d}. In addition, for grid points~$(k^0, \, k^1) \notin \mathcal H$ the entries are scaled by a factor~$-c/4$, where~$c$ is the number of neighbouring grid points in~$\mathcal H$  and~$c = 2 c_1 + c_0$, where $c_0,c_1 \in \{0, \, 1\}$. This is accomplished by the three controlled~$R_Y$ rotation gates with  angles~$\gamma_1 = - 2 \sin^{-1} (1/4)$, ~$\gamma_2 = - 2 \sin^{-1} (1/2)$ and~$\gamma_3 = - 2 \sin^{-1} (3/4)$. In the remaining substeps~$S3$ and~$S4$ all the ancilla qubits are uncomputed. 

The block encoding of the immersed Neumann system matrix~$A_N$ is obtained by subtracting  the block-encoded diagonal correction matrix~$D_{N} (  I \circ A )/ (2 d)$ from the 
 block-encoded immersed system matrix~$A_D$ using the standard LCU construction.

%
\section{Solution of linear systems of equations \label{sec:lsolve}}
%
We next discuss the solution of the discrete system using the quantum singular value transformation (QSVT) framework for computing functions of block-encoded matrices. Here, the function of interest is the inverse of the system matrix. QSVT extends quantum signal processing (QSP), originally developed for univariate polynomial representations, to block-encoded matrices.

%
\subsection{Quantum signal processing (QSP) \label{sec:qsp}}
%
In QSP, a univariate polynomial~$p(x)$, where~$x \in [-1,1]$, is represented as the product of a sequence of unitary  $2 \times 2$ matrices. The QSP sequence consists of unitary matrices of two types, referred to as the \emph{signal operator}  and~\emph{signal-processing operator}. The signal operator encodes the function argument~$x$ and is defined as  
\begin{equation}
	W(x) = \begin{pmatrix}  x & - \sqrt{1-x^2} \\  \sqrt{1-x^2} & x \end{pmatrix} \, .
\end{equation}
This matrix is equivalent to the Pauli~$Y$ rotation matrix   
\begin{equation}
	R_Y (\theta)  =  e^{-i \theta Y /2} = \begin{pmatrix*}[r] \cos  \tfrac{\theta}{2}  & -  \sin \frac{\theta}{2} \\[0.35em]    \sin \tfrac{\theta}{2} & \cos \tfrac{\theta}{2} \end{pmatrix*} \, ,
\end{equation}
with~$\theta  =  2 \cos^{-1} x$. The first key observation is that the $\rho$-th power of~$W(x)$ yields the Chebyshev basis functions of degree~$\rho$. For instance,  
\begin{equation}
	W^{\rho=3} (x) =  \begin{pmatrix}  -3 x + 4 x^3 & (1 - 4 x^2) \sqrt{1-x^2}  \\   (-1 + 4 x^2) \sqrt{1-x^2} & -3 x + 4 x^3   \end{pmatrix}  
		      =  \begin{pmatrix}  C_{T,3}(x) &  - C_{U,2} (x)  \sqrt{1-x^2}\\   C_{U,2}(x)  \sqrt{1-x^2}&  C_{T,3}(x)   \end{pmatrix}  \, . 
\end{equation}
Here,~$C_{T,3}(x)$ is the Chebyshev function of the first kind of degree~$\rho = 3$ and~$C_{U,2}(x)$ is the Chebyshev function of the second kind of degree~$\rho-1 = 2$.

\begin{figure}[!tbp]
    \centering
    \scalebox{1}{\includegraphics{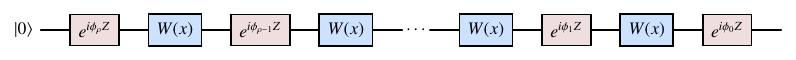}}
    \caption{QSP circuit for a univariate polynomial~$p(x)$ of degree $\rho$. The signal  operator $W(x)$ is fixed and the signal processing operator $e^{i\phi_j Z}$ depends on the phase factors $\phi_j \in \{\phi_0, \, \phi_1, \, \dotsc , \,  \phi_\rho \}$. \label{fig:qsp_circuit}}
\end{figure}

In the QSP construction, polynomials are represented by interleaving the signal operator~$W(x)$ with the signal-processing matrix, which is the Pauli Z rotation matrix
\begin{equation}
	e^{i \phi_j Z} = \begin{pmatrix} e^{i \phi_j } & 0 \\ 0 & e^{-i \phi_j }\end{pmatrix} \, ,
\end{equation}
where~$\phi_j \in \{\phi_0, \, \phi_1, \, \dotsc , \,  \phi_\rho \}$ is a phase factor.  According to Gily\'en et al.~\cite[Theorem 4]{gilyen2019quantum} and Lin~\cite[Theorem 7.20]{lin2022lecture}, the QSP construction reads
\begin{equation}
	U_\phi (x) = e^{i \phi_0 Z} \prod_{j=1}^\rho \left (  W(x) e^{i \phi_j Z} \right ) = \begin{pmatrix}  p(x) & - q(x) \sqrt{1-x^2} \\ q^*(x) \sqrt{1-x^2} & p^*(x)\end{pmatrix} \, , 
\end{equation}
where the star denotes complex conjugation. The two polynomials~$p(x)$ and~$q(x)$ satisfy
\begin{itemize}
	\item[--] $\deg p(x) \le \rho$ and $\deg q(x) \le \rho-1$, 
	\item[--]  $p(x)$ is symmetric and $q(x)$ is antisymmetric when $\rho \in \{2, \, 4, \, \, 6, \, \dotsc \}$, 
	\item[--]  $p(x)$ is antisymmetric and $q(x)$ is symmetric when $\rho \in \{1, \, 3, \, \, 5, \, \dotsc \}$, 
	\item[--] and $| p(x) |^2 + (1-x^2) | q(x)|^2 =1, \quad \forall x \in [-1,1]$. 
\end{itemize}
See Figure~\ref{fig:qsp_circuit} for the circuit implementation of~$U_\phi(x)$. 

\begin{figure}[!tbp]
    \centering
    \includegraphics[scale=0.5]{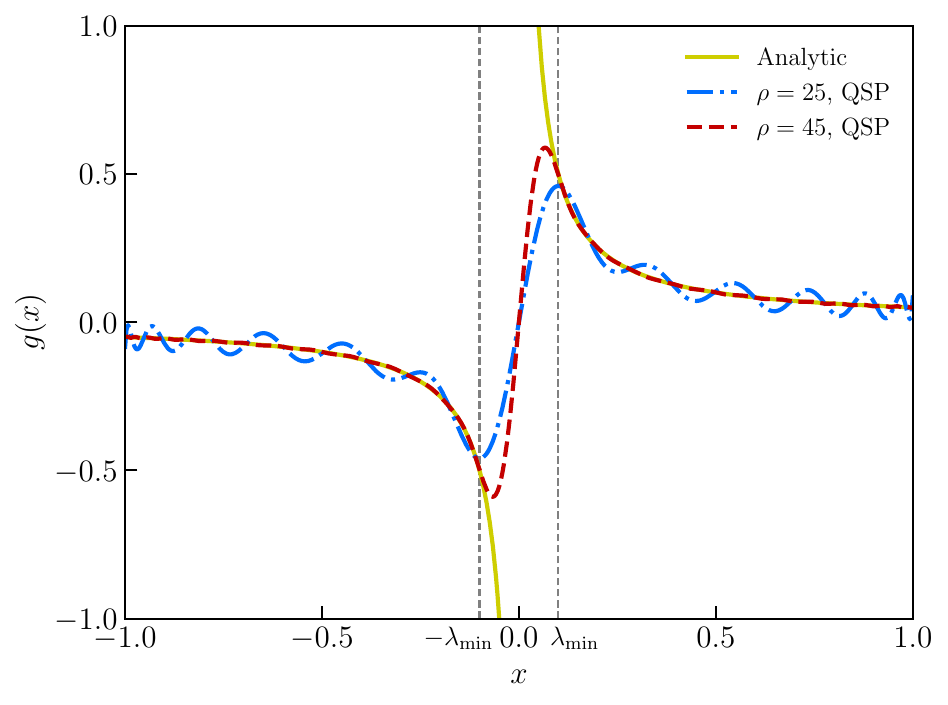}
    \caption{Inverse function~$1/x$ and its scaled Chebyshev polynomial approximations~$\tilde g(x)$ of degrees~$\rho=25$ and~$\rho=45$ for~$\lambda_{\text{min}}=0.1$. The polynomials are obtained by fitting to the regularised inverse function~$g(x)=\lambda_{\text{min}}/(2 x)$ on $S_{\! \lambda}  = [-1, \, -\lambda_{\text{min}}]\cup[\lambda_{\text{min}}, \, 1]$.\label{fig:qsp-approx-plot}}
\end{figure}

The QSP construction makes it possible to approximately quantum encode any function by approximating it with a polynomial and determining the phase factors~$\{ \phi_j \}_{j=0}^\rho$ on a classical computer. 
In solving the linear system of equations, the QSP approximation of the inverse function $g(x)=1/x$ is required. Due to the singularity of~$g(x)$ at~$x =0 $, only the interval $S_{\! \lambda} = [-1, \, -\lambda_{\text{min}}] \cup [\lambda_{\text{min}}, \, 1]$ is considered. The reason for referring to the interval boundary as~$\lambda_{\min}$ will become clear below.  Because of the requirement~$|g(x) | < 1$, we consider on~$S_{\! \lambda}$ the rescaled inverse function
\begin{equation}\label{eq:qsp-approx-func}
	g(x) = \frac{\lambda_{\text{min}}}{2 x} \, . 
\end{equation}
We determine the phase factors~$\{ \phi_j \}_{j=0}^\rho$ using the QSPPACK library~\cite{dong2021qsppack}. In QSPPACK, $g(x)$ is first approximated by a Chebyshev polynomial~$\tilde g(x)$ of a given degree~$\rho$ using standard approximation techniques. The Chebyshev approximation error~$\epsilon$ is defined as 
\begin{equation}
	\epsilon = \max_{S_{\! \lambda}} | g(x) - \tilde g(x) | \, . 
\end{equation}
According to the QSP construction,~$\rho$ must be odd because~$g(x)$ and its approximation~$\tilde g(x)$ are antisymmetric. The obtained Chebyshev polynomial~$\tilde g(x)$ is exactly representable by QSP.  The corresponding phase factors~$\{ \phi_j \}_{j=0}^\rho$ are determined by minimising the error in a chosen norm, i.e.,  
\begin{equation}
	\{ \phi_j^* \}  = \argmin_{ \{ \phi_j \} } \left \| \operatorname{Re} \left ( \bra{0}  U_{\phi} (x)\ket 0 \right  ) - \tilde g(x) \right  \|  \, , 
\end{equation}
where $\operatorname{Re}$ denotes the real part of a complex number. In Figure~\ref{fig:qsp-approx-plot}, the QSP approximation of~$g(x)$ with polynomials of degree~$\rho=25$ and~$\rho=45$ is depicted. As expected, an increase in the polynomial degree~$\rho$ leads to a reduction in the approximation error. 

%
\subsection{Quantum singular value transformation (QSVT) \label{sec:qsvt}}
%
We consider the eigendecomposition of the finite element system matrix and its action on specific eigenvectors to illustrate the QSVT construction. All the system matrices introduced in this paper are real-valued and symmetric. Those associated with simply connected domains are positive definite, whereas those associated with immersed domains are only positive semidefinite because of the presence of inactive grid points outside the domain, leading to zero eigenvalues. 

The system matrices admit the eigendecomposition 
\begin{equation}
	A/\sigma  = \sum_j \lambda_j \ket {w_j} \bra {w_j} \, , 
\end{equation}
where~$\lambda_j$ and~$w_j$ are corresponding eigenvalues and eigenvectors, and~$\sigma$ is the subnormalisation factor of the block-encoding.  The minimum non-zero and maximum eigenvalues are denoted by~$\lambda_{\text{min}}$ and~$\lambda_{\text{max}}$, respectively, and the maximum eigenvalue satisfies~$\lambda_{\text{max}} \le 1$. The scaled inverse is approximated using the Chebyshev polynomial~$\tilde g(\lambda)$ approximating the rescaled inverse function~$g(\lambda)$ defined in~\eqref{eq:qsp-approx-func}, so that
\begin{equation}
	\frac{\lambda_{\text{min}} \sigma}{2} A^{-1}  = \sum_j \frac{\lambda_{\text{min}}}{2  \lambda_j}\ket {w_j} \bra {w_j}  \approx  \sum_j \tilde g (\lambda_j) \ket {w_j} \bra {w_j} \, .
\end{equation}
Here, \(A^{-1}\) denotes the inverse restricted to active grid points inside the domain, and both sums are only over the non-zero eigenvalues. The minimum non-zero eigenvalue of~$A/\sigma$ required to determine~$\tilde g(\lambda)$ can be classically estimated using standard techniques.

Each system matrix introduced in this paper has a corresponding block-encoding~$U_A$ given by the associated quantum circuit. We write for $U_A$ and its adjoint~$U_A^\dagger$ in block-matrix form
\begin{equation}
	U_{A} = \begin{pmatrix} A/\sigma & * \\ * & * \end{pmatrix} \, , \qquad U_{A}^\dagger = \begin{pmatrix} A/\sigma & *^\dagger \\ *^\dagger & *^\dagger \end{pmatrix} \, .
\end{equation}
The overall dimension of~$U_A$ and~$U_A^\dagger$ depends on the additional qubits beyond the ones in the registers~$\ket{k^0}$ in 1D and~$\ket{k^0} \ket{k^1}$ in 2D used for labelling the grid points. The matrix-vector products~$U_{A} \ket 0 \ket{w_j}$  and~$U_{A}^\dagger \ket 0 \ket{w_j}$ expressed as block matrix operations read 
\begin{equation}
	U_A 
	\begin{pmatrix}
		\ket{w_j} \\ 0 
	\end{pmatrix} 
	= \begin{pmatrix} 
		\lambda_j \ket{w_j} \\ 
		\sqrt{1-\lambda_j^2} \ket{w_{\perp, j}} 
	   \end{pmatrix} \, ,  \qquad  	
	U_A^\dagger
	\begin{pmatrix}
		\ket{w_j} \\ 0 
	\end{pmatrix} 
	= \begin{pmatrix} 
		\lambda_j \ket{w_j} \\ 
		\sqrt{1-\lambda_j^2} \ket{w'_{\perp, j}} 
	   \end{pmatrix} \, , 
\end{equation} 
where~$\ket{w_{\perp, j}}$ and~$\ket{w'_{\perp, j}}$ are  two unit vectors.  The application of~$U_A^\dagger$ on the first equation and~$U_A$ on the second equation, noting that~$U_A U_A^\dagger=I$, yields
\begin{equation}
	U_A^\dagger \begin{pmatrix}  0 \\ \ket{w_{\perp, j}}  \end{pmatrix} =  \begin{pmatrix}   \sqrt{1 - \lambda_j^2}  \ket {w_{j}} \\ -\lambda_j \ket{w'_{\perp,j}} \end{pmatrix} \, , \qquad 
	U_A \begin{pmatrix}  0 \\ \ket{w'_{\perp, j}}  \end{pmatrix} =  \begin{pmatrix}   \sqrt{1 - \lambda_j^2}  \ket {w_{j}} \\ -\lambda_j \ket{w_{\perp,j}} \end{pmatrix} \, . 
\end{equation}
Hence,~$U_A$ and~$U_A^\dagger$ act in the two-dimensional subspaces defined by~$\{\ket {w_j}, \, \ket{w_{\perp, j}}  \}$  and ~$\{\ket {w_j}, \, \ket{w'_{\perp, j}}  \}$, and have the following entries
\begin{equation} \label{eq:UandUdagger}
	U_A 
	\begin{pmatrix}
		\ket{w_j} \\  \ket{w'_{\perp, j}}  
	\end{pmatrix}  = 
	\begin{pmatrix}
		\lambda_j & \sqrt{1- \lambda_j^2}  \\
		 \sqrt{1- \lambda_j^2} & - \lambda_j 
	\end{pmatrix}
		\begin{pmatrix}
		\ket{w_j} \\  \ket{w_{\perp, j}}  
	\end{pmatrix}  \, , \quad 
		U_A ^\dagger
	\begin{pmatrix}
		\ket{w_j} \\  \ket{w_{\perp, j}}  
	\end{pmatrix}  = 
	\begin{pmatrix}
		\lambda_j & \sqrt{1- \lambda_j^2}  \\
		 \sqrt{1- \lambda_j^2} & - \lambda_j 
	\end{pmatrix}
		\begin{pmatrix}
		\ket{w_j} \\  \ket{w'_{\perp, j}}  
	\end{pmatrix}  \, . 
\end{equation}
Next, to make use of the QSP construction,  a composite unitary is defined, which acts in the subspace~$\{ \ket{w_j}, \ket{w'_{\perp, j}} \}$, and yields a~$2\times 2$ QSP signal operator.  By making use of~\eqref{eq:UandUdagger} we can write
\begin{equation}
	U_A^\dagger Z_{\Pi}  U_A Z_{\Pi}   \begin{pmatrix}
		\ket{w_j} \\  \ket{w'_{\perp, j}}  
	\end{pmatrix} = 
	\begin{pmatrix}
		\lambda_j & - \sqrt{1- \lambda_j^2}  \\
		 \sqrt{1- \lambda_j^2} &  \lambda_j 
	\end{pmatrix}^2	
	 \begin{pmatrix}
		\ket{w_j} \\  \ket{w'_{\perp, j}}  
	\end{pmatrix} \, ,
\end{equation}
where~$Z_{\Pi}$ is the reflection matrix  
\begin{equation}
	Z_{\Pi} = \begin{pmatrix*}[r] I & 0 \\ 0 & -I\end{pmatrix*} \, . 
\end{equation}
According to the QSP construction, in the considered subspace, a Chebyshev function of even degree~$\rho \ge 2$ can be obtained as follows
\begin{equation}
	(U_A^\dagger Z_{\Pi}  U_A Z_{\Pi})^{\rho/2}  = 
	\begin{pmatrix}
		\lambda_j & - \sqrt{1- \lambda_j^2}  \\
		 \sqrt{1- \lambda_j^2} &  \lambda_j 
	\end{pmatrix}^\rho	 =
	 \begin{pmatrix}  C_{T,\rho}(\lambda_j) &  - C_{U,\rho-1} (\lambda_j)  \sqrt{1-\lambda_j^2}\\   C_{U,\rho-1}(\lambda_j)  \sqrt{1-\lambda_j^2}&  C_{T,\rho}(\lambda_j)   \end{pmatrix} \, .
\end{equation}
Chebyshev functions of odd degree require a slightly different composite unitary~\cite{gilyen2019quantum,lin2022lecture}.  The sketched construction allows for computing the action of the Chebyshev matrix function~$C_{T, \rho} (A/\sigma) $ on any vector~$\ket {f}$, i.e. $C_{T, \rho} (A/\sigma)  \ket f$ by evaluating~$(U_A^\dagger Z_{\Pi}  U_A Z_{\Pi})^{\rho/2} \ket 0 \ket f $.   

As in QSP, arbitrary polynomials are obtained by interleaving the signal matrix with the signal-processing matrix. Focusing on antisymmetric polynomials, a polynomial of odd degree~$\rho$ with phase factors~$\phi = \{ \phi_0, \, \phi_1, \, \dotsc , \, \phi_{\rho}  \} $, the QSVT construction reads
\begin{equation}	
	U_\phi (A/\sigma ) = (-i)^\rho e^{i \phi_0 Z_{\Pi}} U_A e^{i \phi_1 Z_{\Pi}} \prod_{j=1}^{(\rho-1)/2} \left ( U_A^\dagger  e^{i \phi_{2j} Z_{\Pi}} U_A  e^{i \phi_{2j+1} Z_\Pi}   \right )  = \begin{pmatrix} p(A/\sigma) & * \\ * & * \end{pmatrix}\, .
\end{equation}
\begin{figure*}[!tbp]
    \centering
    \scalebox{1}{\includegraphics{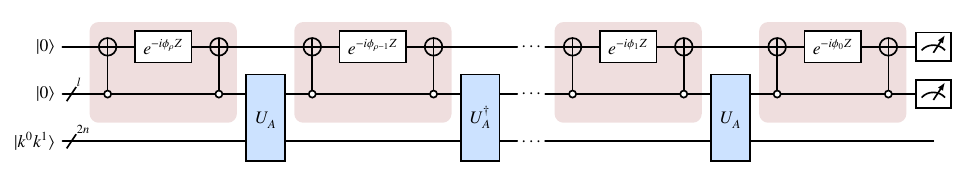}}
    \caption{QSVT circuit of the matrix polynomial~$p(A)$ of degree~$\rho$ of the block-encoded matrix~$A$.  The block-encoding unitary $U_A$ and its conjugate transpose $U_A^\dagger$ are applied alternately. Between the block-encoding unitaries, a projector-controlled Pauli Z rotation gate is applied to the ancilla qubits $\ket{0}\ket{0}^{\otimes l}$. 
    \label{fig:qsvt_circuit}}
\end{figure*}
See Figure~\ref{fig:qsvt_circuit} for the circuit implementation of~$U_\phi(A/\sigma)$. 

%
\section{Examples \label{sec:examples}}
We present a series of numerical experiments to assess the computational complexity and verify the accuracy of the proposed {\em IQFEM} framework. All the introduced circuits are implemented in Qiskit and executed using a noiseless statevector simulator~\cite{qiskit2024}. To assess computational complexity, we let Qiskit decompose the circuits into two-qubit $\mathit{CNOT}$ and single-qubit rotation $U_3$ gates. The gate set $\{ \mathit{CNOT}, U_3\}$ is universal, meaning that any unitary operation can be implemented using only those two types of gates. We first study two one-dimensional examples with known analytical solutions: a heterogeneous problem and an immersed problem with Dirichlet boundary conditions.  Subsequently, we consider three two-dimensional examples: a heterogeneous problem and an immersed problem with Dirichlet and Neumann boundary conditions. In all the examples, the system equations are solved using the QSVT-based quantum linear system solver, and the required phase factors are determined using the QSPPACK library~\cite{dong2021qsppack}. 
%
\subsection{One-dimensional problems}
%
\subsubsection{Heterogeneous diffusion}
\begin{figure}[!tbp]
    \centering
    \subfloat[Source~$f(x)$\label{fig:1d-het-figure_a}]{%
        \begin{minipage}[b]{0.43\textwidth}
            \centering
            \includegraphics[width=\linewidth, keepaspectratio]{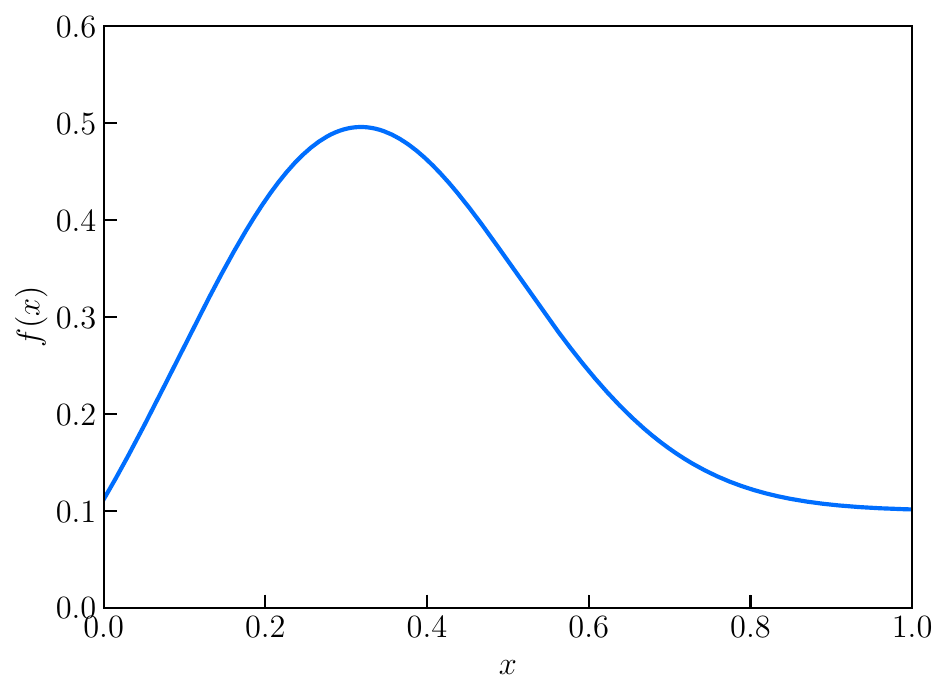}
        \end{minipage}
    }
   \hspace{0.075\textwidth}
    \subfloat[Solution~$u(x)$\label{fig:1d-het-figure_b}]{%
        \begin{minipage}[b]{0.43\textwidth}
            \centering
            \includegraphics[width=\linewidth, keepaspectratio]{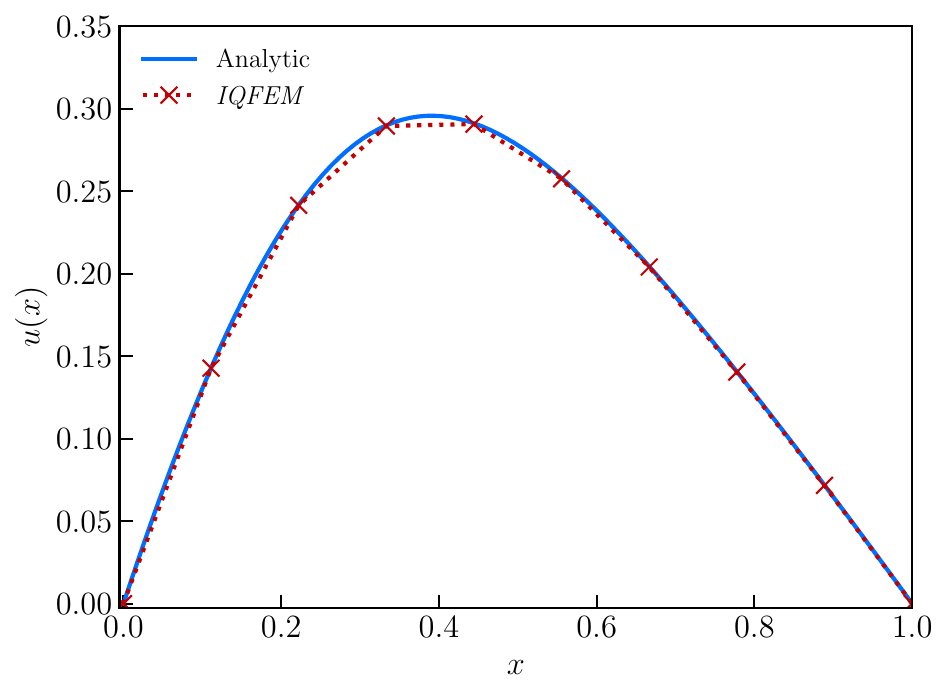}
        \end{minipage}
    }
    \vfill
    \subfloat[Number of~$U_3$ and~$\mathit{CNOT}$ gates\label{fig:1d-het-figure_d}]{%
        \begin{minipage}[b]{0.43\textwidth}
            \centering
            \includegraphics[width=1.0\linewidth, keepaspectratio]{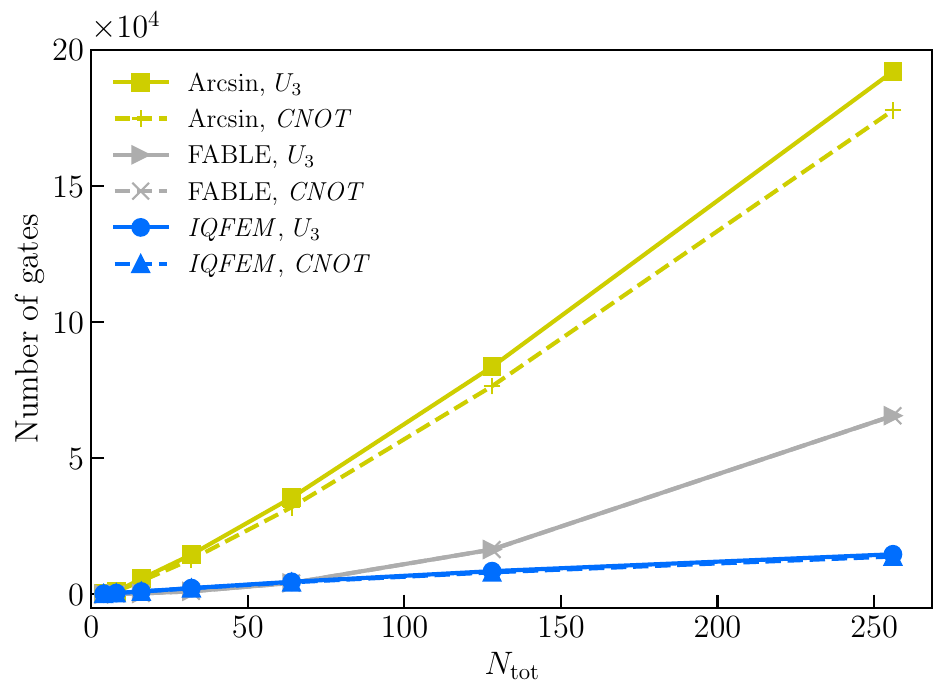}
        \end{minipage}
    }
    \hspace{0.075\textwidth}
    \subfloat[Relative $L^2$-norm error in $u_h(x)$\label{fig:1d-het-figure_c}]{%
        \begin{minipage}[b]{0.43\textwidth}
            \centering
            \includegraphics[width=0.98\linewidth, keepaspectratio]{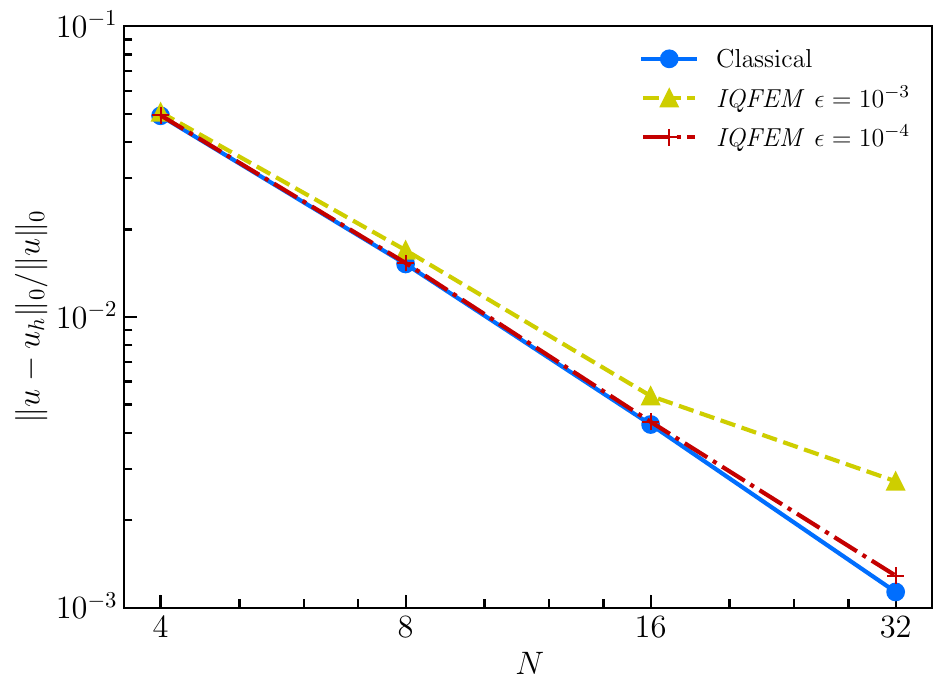}
        \end{minipage}
    }
      \caption{One-dimensional heterogeneous Poisson problem. (a) Source term. (b) Analytical and {\em IQFEM} solutions for $N = 2^3$ grid points. (c) Number of required~$U_3$ and~$\mathit{CNOT}$ gates to block-encode the system matrix using the {\em IQFEM}, FABLE and Arcsin approaches. (d) Convergence of the relative $L^2$-norm error of the {\em IQFEM} solution for QSVT approximation errors $\epsilon \in \{10^{-3}, \, 10^{-4}\}$. }
    \label{fig:1d-het-figure}
\end{figure}

%
We consider on the domain~$\Omega=(0,1)$ the heterogeneous Poisson problem subject to the Dirichlet boundary conditions~$u(0)=u(1)=0$. The spatially varying diffusion coefficient is given by 
\begin{equation}
    \mu(x)=\frac{1+x}{10} \ .
\end{equation}
The source term $f(x)$ is  chosen as 
\begin{equation}
    f(x)=0.1 \operatorname{erf}\!\left(\frac{20x-5}{6}\right) + \frac{2(1+x)}{3\sqrt{\pi}} \exp\!\left[-\left(\frac{20x-5}{6}\right)^2\right] \, .
\end{equation}
The corresponding analytical solution is given by 
\begin{equation}
    u(x)= c_0 + c_1 \ln\!\left(\frac{1+x}{10}\right) - (x-0.25)\operatorname{erf}\!\left(\frac{20x-5}{6}\right) - \frac{3}{10\sqrt{\pi}}\exp\!\left[-\left(\frac{20x-5}{6}\right)^2\right],
\end{equation}
where the integration constants $c_0$ and $c_1$ are chosen to satisfy the Dirichlet boundary conditions. The source term and the analytical and {\em IQFEM} solutions for~$N=2^3$ internal grid points are plotted in Figures~\ref{fig:1d-het-figure_a} and~\ref{fig:1d-het-figure_b}, respectively.

The~$U_3$ and~$\mathit{CNOT}$ gate counts required to block-encode the system matrix as a function of the number of grid points~$N$ are plotted in Figure~\ref{fig:1d-het-figure_d}. In {\em IQFEM}, the number of gates scales polylogarithmically with the number of grid points. Figure~\ref{fig:1d-het-figure_d} also includes the gate counts for the Arcsin and FABLE block encoding methods \cite{lapworth2024encoding,camps2022fable}, which exhibit a polynomial scaling. It is worth noting that both methods require a significant number of gates even for relatively modest grid sizes, far exceeding the capabilities anticipated for near-term quantum computing platforms. 

The system is solved using the introduced QSVT-based quantum linear system solver.  The tolerance controlling the degree of the Chebyshev approximation to the regularised inverse function \mbox{$g(x)= \lambda_\text{min}/(2 x)$} is chosen as either \mbox{$\epsilon=10^{-3}$} or \mbox{$\epsilon=10^{-4}$}.  The convergence of the relative $L^2$-norm error of the quantum solution for the two different values of~$\epsilon$ is compared with that of the classical finite element solution in Figure~\ref{fig:1d-het-figure_c}. Choosing a sufficiently small value of~$\epsilon$ for a given grid resolution recovers the optimal convergence rate.  

Furthermore, {\em IQFEM} has a subnormalisation factor $\sigma = 3 \| A \|_{\text{max}} $, see~\eqref{eq:Av_blockMatrixForm}, in contrast to the subnormalisation factor~$\sigma= N \| A \|_{\text{max}}$ of the FABLE and Arcsin encodings. A larger subnormalisation factor reduces the smallest non-zero eigenvalue~$\lambda_{\text{min}}$, which leads to an increase in the polynomial degree and the number of gates in the QSVT matrix inversion.

%
\subsubsection{Immersed Dirichlet boundary}
\begin{figure}[!tbp]
    \centering
    \subfloat[Source~$f(x)$\label{fig:1d-im-figure_a}]{%
        \begin{minipage}[b]{0.43\textwidth}
            \centering
            \includegraphics[width=\linewidth, keepaspectratio]{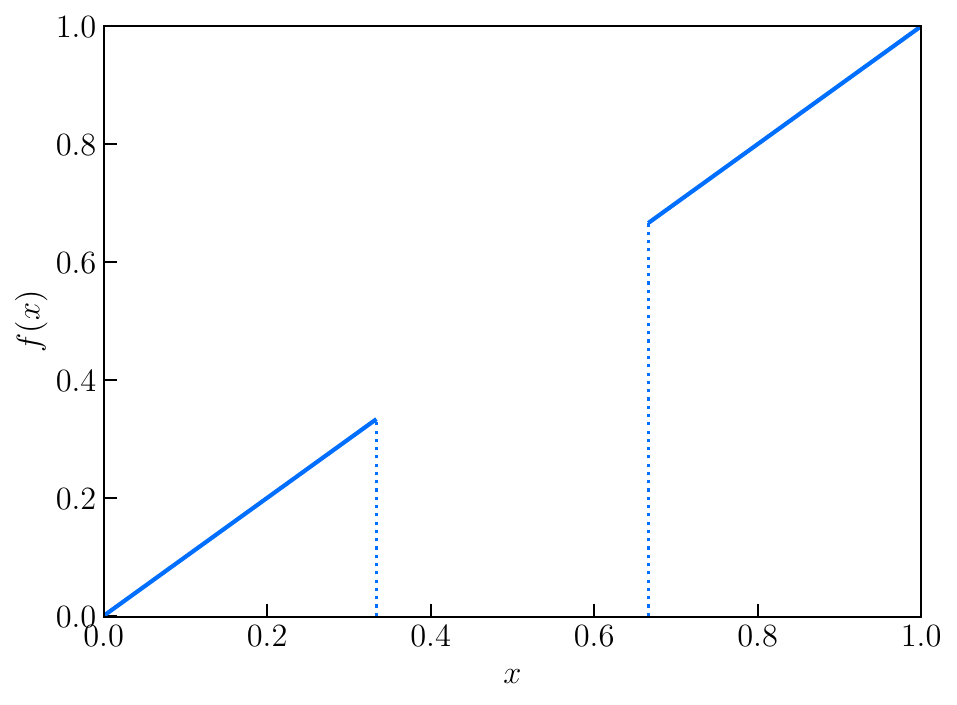}
        \end{minipage}
    }
    \hspace{0.075\textwidth}
    \subfloat[Solution~$u(x)$\label{fig:1d-im-figure_b}]{%
        \begin{minipage}[b]{0.43\textwidth}
            \centering
            \includegraphics[width=\linewidth, keepaspectratio]{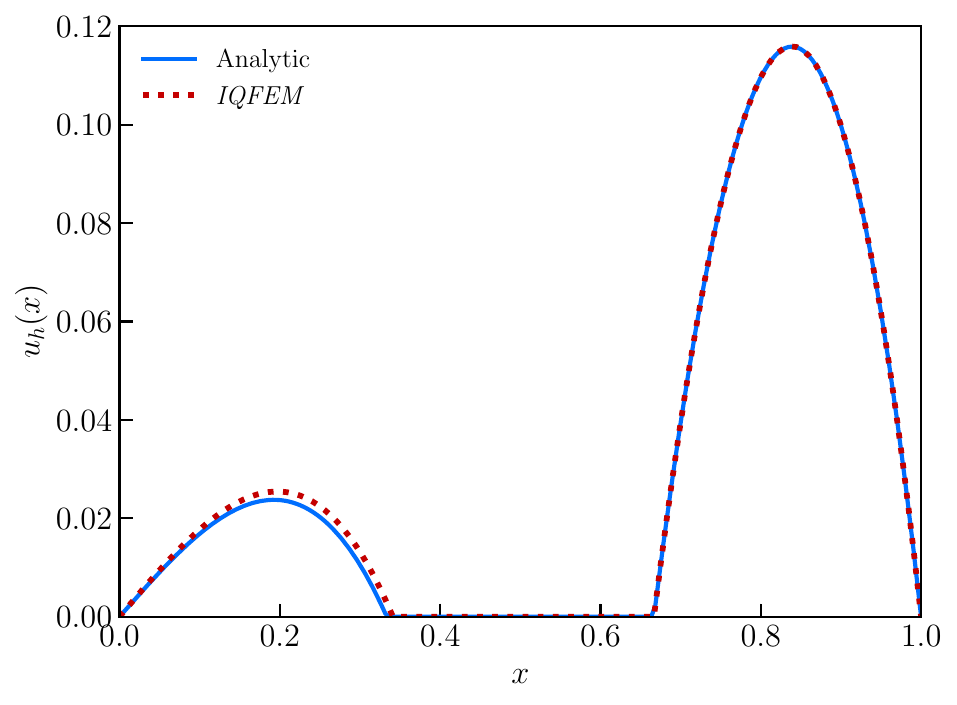}
        \end{minipage}
    }
    \vspace{0.2cm}
    \subfloat[Number of~$U_3$ and~$\mathit{CNOT}$ gates\label{fig:1d-im-figure_d}]{%
        \begin{minipage}[b]{0.43\textwidth}
            \centering
            \includegraphics[width=\linewidth, keepaspectratio]{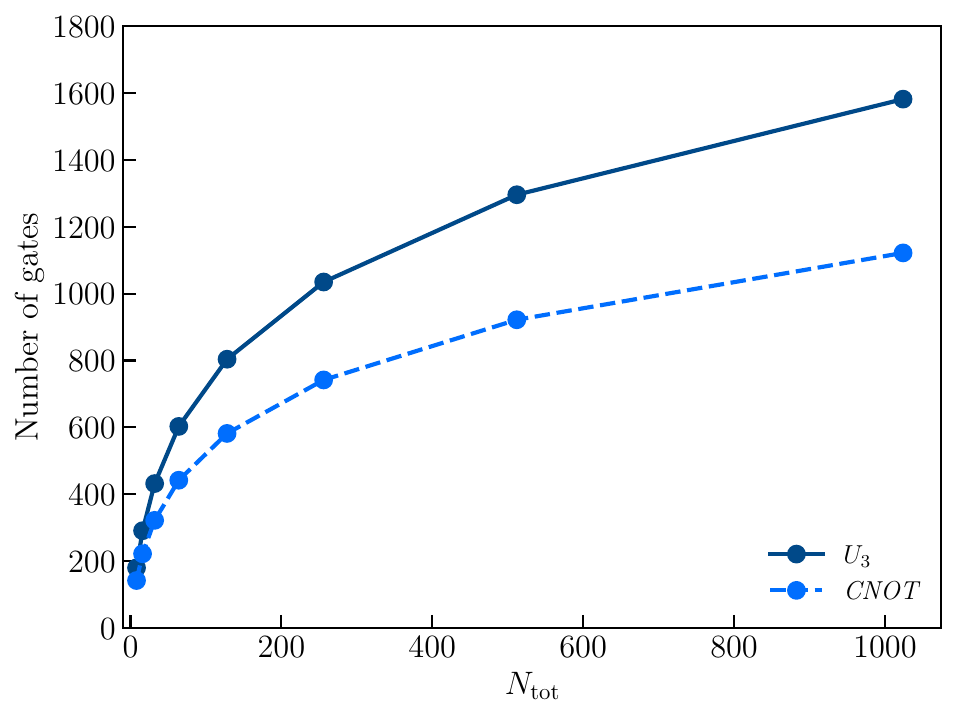}
        \end{minipage}
    }
    \hspace{0.075\textwidth}
    \subfloat[Relative $L^2$-norm error in $u_h(x)$\label{fig:1d-im-figure_c}]{%
        \begin{minipage}[b]{0.43\textwidth}
            \centering
            \includegraphics[width=\linewidth, keepaspectratio]{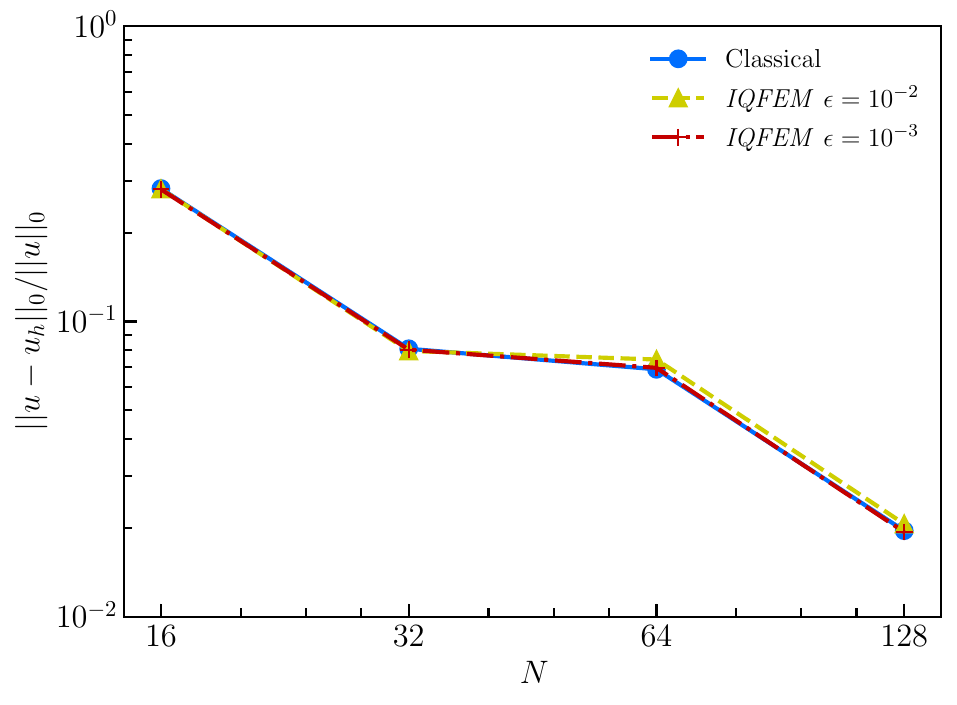}
        \end{minipage}
    }
    \hfill
        \caption{One-dimensional Poisson problem with an immersed hole. (a) Source term. (b) Analytical and {\em IQFEM} solutions for $N = 2^7$ grid points. (c) Number of required~$U_3$ and~$\mathit{CNOT}$ gates to block-encode the system matrix. (d) Convergence of the relative $L^2$-norm error of the {\em IQFEM} solution for QSVT approximation errors $\epsilon \in \{10^{-2}, 10^{-3}\}$.}
    \label{fig:1d-im-figure}
\end{figure}

%
We now consider a Poisson problem defined on the disconnected domain~ $\Omega=(0, \, 1)\setminus [1/3, \, 2/3]$. The boundary conditions at the external and the immersed boundaries are~$u(0)=u(1)=0$ and~$u(1/3)=u(2/3)=0$, respectively. The diffusion coefficient is constant~$\mu(x)=0.1$, and the forcing is linear $f(x)=x$. The corresponding analytical solution is given by
\begin{equation}
    u(x)=
    \begin{cases}
        \dfrac{10x-90x^3}{54}\, ,       & x \in (0, \, 1/3) \, ,   \\[0.5em]
        \dfrac{-100+190x-90x^3}{54}\, , & x \in (2/3, \, 1)\, .
    \end{cases}
\end{equation}
The source term and the analytical and {\em IQFEM} solutions for~$N=2^7$ internal grid points are plotted in Figures~\ref{fig:1d-im-figure_a} and~\ref{fig:1d-im-figure_b}, respectively. 

In {\em IQFEM}, the disconnected domain is represented as a domain with an internal hole. The geometry of the hole is defined by a level-set function introduced in Section~\ref{sec:oneDlsetAndIndic}. The domain~$(0, \,1)$ is discretised into~$N+1$ equal intervals. The Dirichlet boundary conditions~$u(0)=u(1)=0$ are enforced exactly, while Dirichlet boundary conditions on the immersed boundary~$u(1/3)=u(2/3)=0$  are enforced by eliminating the system matrix entries corresponding to the grid points in the hole. This is accomplished using the indicator function, which is evaluated using quantum arithmetic. 

The {\em IQFEM} grids are not nested because increasing the number of qubits~$n$ to represent the grid points by one doubles the number of grid points. However, a uniform refinement of the grid with~$N$ interior points yields a grid with~$2N+1$ interior points. Consequently, the immersed boundary location is not fixed between successive {\em IQFEM} grids, which leads to approximation errors because boundary conditions are enforced only at grid points. 

The~$U_3$ and $\mathit{CNOT}$ gate counts required to block-encode the immersed system matrix are plotted in Figure~\ref{fig:1d-im-figure_d}. As in the previous example, they are relatively small and exhibit a polylogarithmic scaling. These counts include all gates required to evaluate the level-set and the indicator functions.  The convergence of the $L^2$-norm error in {\em IQFEM} solutions obtained with classical and QSVT-based linear equation solvers and tolerances~$\epsilon=10^{-2}$ and~$\epsilon=10^{-3}$ is plotted in Figure~\ref{fig:1d-im-figure_c}.  The $L^2$-norm exhibits approximately linear convergence, which is suboptimal due to changes in the approximated boundary location under grid refinement. 

\begin{figure}[!tbp]
    \centering
    \includegraphics[scale=0.43]{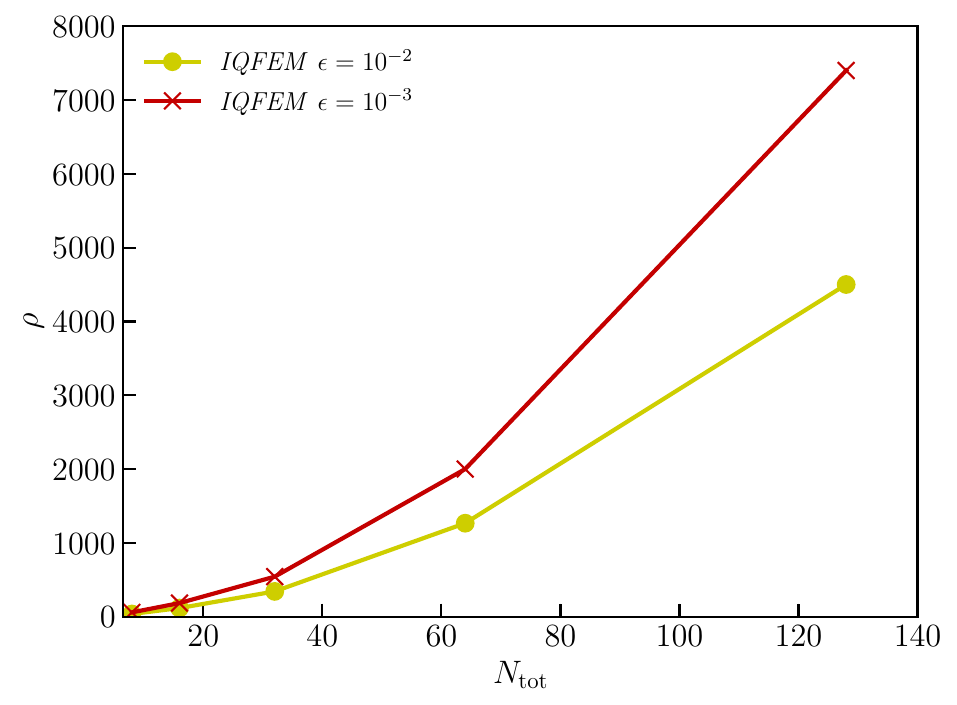}
    \vspace{-0.2cm}
    \caption{One-dimensional Poisson problem with immersed Dirichlet boundary. Polynomial degree $\rho$ versus number of grid points for different QSVT approximation errors $\epsilon \in \{10^{-2}, 10^{-3}\}$.\label{fig:1d-im-figure_poly-degree}}
\end{figure}

The degree~$\rho$ of the Chebyshev polynomial approximating the rescaled inverse function in the QSVT solver is shown in Figure~\ref{fig:1d-im-figure_poly-degree}. For a given discretisation with $N$ grid points, we choose $\rho$ such that the maximum approximation error is either~$\epsilon=10^{-2}$ or~$\epsilon=10^{-3}$. The minimum non-zero eigenvalue of the scaled system matrix scales as~$N^{-2}$; therefore, grid refinement leads to an approximately quadratic increase in the required polynomial degree. The Chebyshev approximation error and the finite element discretisation error must be balanced to retain the optimal convergence of the finite element solution. Furthermore, the polynomial degree~$\rho$ is about equal to the number of queries to the block-encoded matrix. Hence, the QSVT gate counts are proportional to~$\rho$ times the gate counts for encoding a single system matrix.   

%
\subsection{Two-dimensional problems}
%
\subsubsection{Heterogeneous diffusion}
\begin{figure}[!tbp]
    \centering
   \subfloat[Diffusion coefficient~$\mu(x_0, x_1)$\label{fig:2d-het-figure_b}]{%
        \begin{minipage}[b]{0.43\textwidth}
            \includegraphics[width=\linewidth, keepaspectratio]{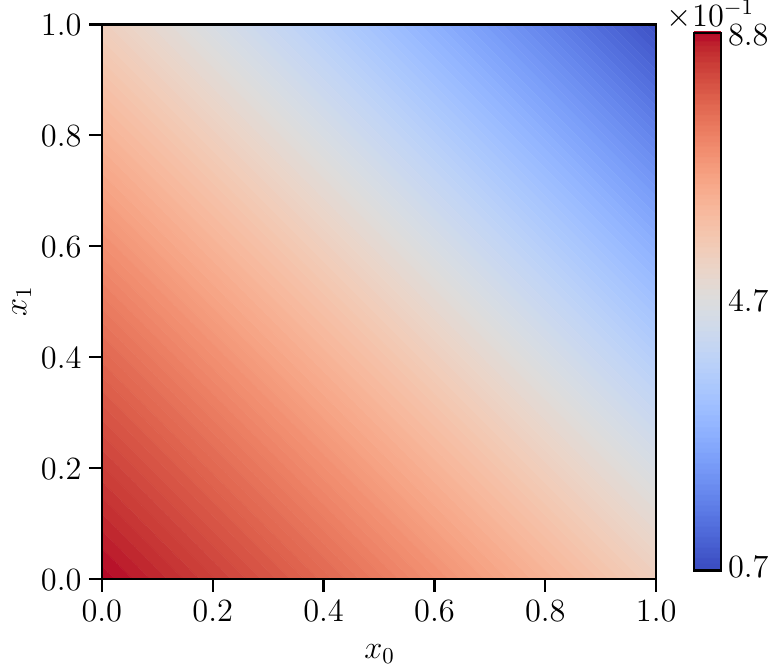}
        \end{minipage}
   }
\hspace{0.075\textwidth}
    \subfloat[{\em IQFEM} solution~$u_h(x_0, x_1)$\label{fig:2d-het-figure_c}]{%
        \begin{minipage}[b]{0.43\textwidth}
            \includegraphics[width=1.035\linewidth, keepaspectratio]{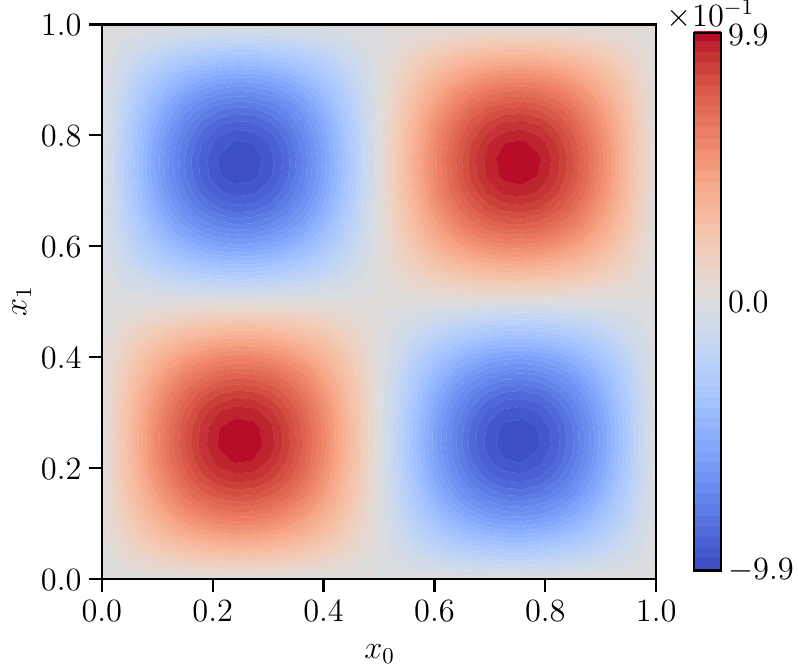}
        \end{minipage}
    }
    \hspace{0.075\textwidth}
    \subfloat[Number of~$U_3$ and~$\mathit{CNOT}$ gates\label{fig:2d-het-figure_f}]{%
        \begin{minipage}[b]{0.43\textwidth}
            \centering
            \includegraphics[width=\linewidth, keepaspectratio]{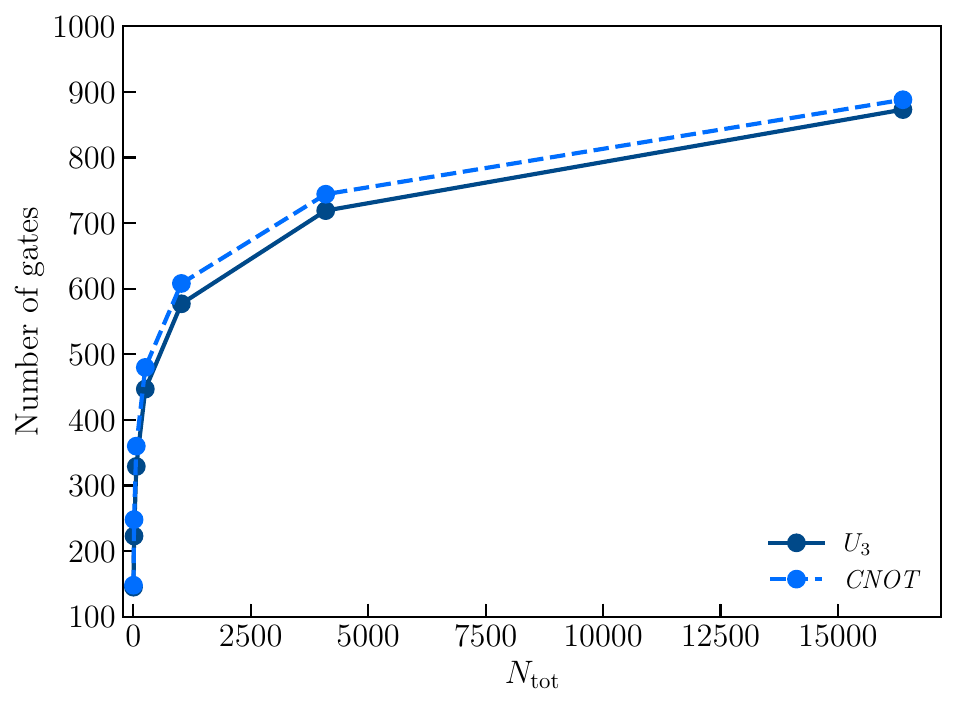}
        \end{minipage}
    }
    \hspace{0.075\textwidth}
    \subfloat[Relative $L^2$-norm error in~$u_h(x_0, x_1)$\label{fig:2d-het-figure_e}]{%
        \begin{minipage}[b]{0.43\textwidth}
            \centering
            \includegraphics[width=\linewidth, keepaspectratio]{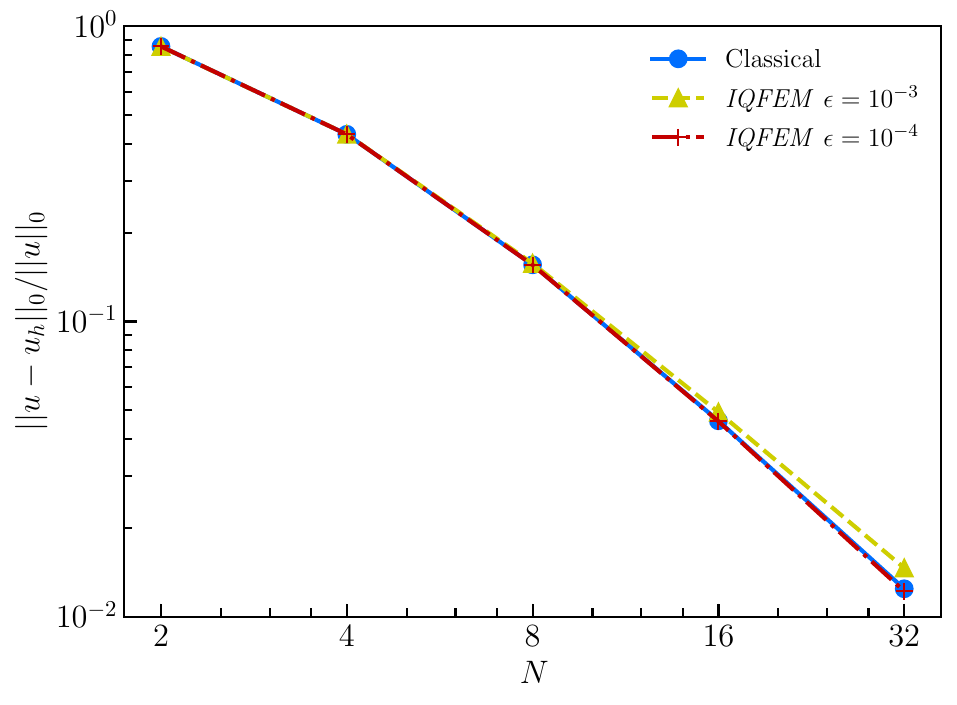}
        \end{minipage}
    }
    \caption{Two-dimensional heterogeneous Poisson problem.  (a) Spatially varying diffusion coefficient~$\mu(x_0, x_1)$. (b) {\em IQFEM} solution~$u_h(x_0, x_1)$ for $2^5 \times 2^5$ grid points. (c) Number of required~$U_3$ and~$\mathit{CNOT}$ gates to block-encode the system matrix. (d) Convergence of the relative $L^2$-norm error of the {\em IQFEM} solution for QSVT approximation errors $\epsilon \in \{10^{-3}, 10^{-4}\}$. }
    \label{fig:2d-het-figure_2}
\end{figure}

%
We proceed to two-dimensional problems and consider the Poisson equation on the unit square $ (0,1)^2$ with homogeneous Dirichlet boundary conditions.
The spatially varying diffusion coefficient is chosen as
\begin{equation}
    \mu(x_0, x_1) =
    \cos\left(\frac{1 + x_0 + x_1}{2}\right) \, , 
\end{equation}
and the source term as
\begin{equation}
    f(x_0, x_1) =
    \pi \sin\left(\frac{1 + x_0 + x_1}{2}\right)
    \sin\left(2\pi(x_0 + x_1)\right)
    + 8\pi^2
    \cos\left(\frac{1 + x_0 + x_1}{2}\right)
    \sin(2\pi x_0)\sin(2\pi x_1) \ .
\end{equation}
The analytical solution to this problem is given by
\begin{equation}
    u(x_0, x_1) = \sin(2 \pi x_0)\sin(2 \pi x_1) \ .
\end{equation}
The isocontours of the spatially varying diffusion coefficient and the {\em IQFEM} solutions on a grid with~$2^5 \times 2^5$ grid points are shown in Figures~\ref{fig:2d-het-figure_b} and \ref{fig:2d-het-figure_c}, respectively. 

The~$U_3$ and~$\mathit{CNOT}$ gate counts required to block-encode the~$N_{\text{tot}}= N \times N$ system matrix are shown in Figure~\ref{fig:2d-het-figure_f}. The total number of required gates depends on the degree of the polynomial used to approximate the diffusion coefficient. The given harmonic diffusion coefficient~$\mu(x_0,x_1)$ can be directly encoded using a sequence of~$R_Y$ gates, see e.g.~\cite{liu2024towards}. The gate-count scaling in {\em IQFEM} is always polylogarithmic, independent of the polynomial degree for approximation.  The convergence of the $L^2$-norm errors in {\em IQFEM} solutions obtained with classical and QSVT-based linear equation solvers and tolerances~$\epsilon=10^{-3}$ and~$\epsilon=10^{-4}$ is shown in Figure~\ref{fig:2d-het-figure_e}.  The QSVT-based solution recovers the optimal quadratic convergence rate when the tolerance is sufficiently small.

%
\subsubsection{Immersed Dirichlet boundary}
%
We next consider a multiply connected domain with an immersed diamond-shaped hole. The hole is centred at $x_0=1/2$ and~$x_1 =1/2$ and has half-width~$1/4$. Its geometry within the domain~$(0, \,1) \times (0, \, 1)$ is described by a level-set function. The diffusion coefficient~$\mu(x_0,x_1)=1$ and the source~$f(x_0,x_1)=1$ are constant.  Zero Dirichlet boundary conditions are imposed on both the external boundary and the hole boundary.

\begin{figure}[!tbp]
    \centering
    \subfloat[{\em IQFEM} solution~$u_h(x_0, x_1)$\label{fig:2d-im-figure_a}]{%
        \begin{minipage}[c]{0.43\textwidth}
            \includegraphics[width=\linewidth, keepaspectratio]{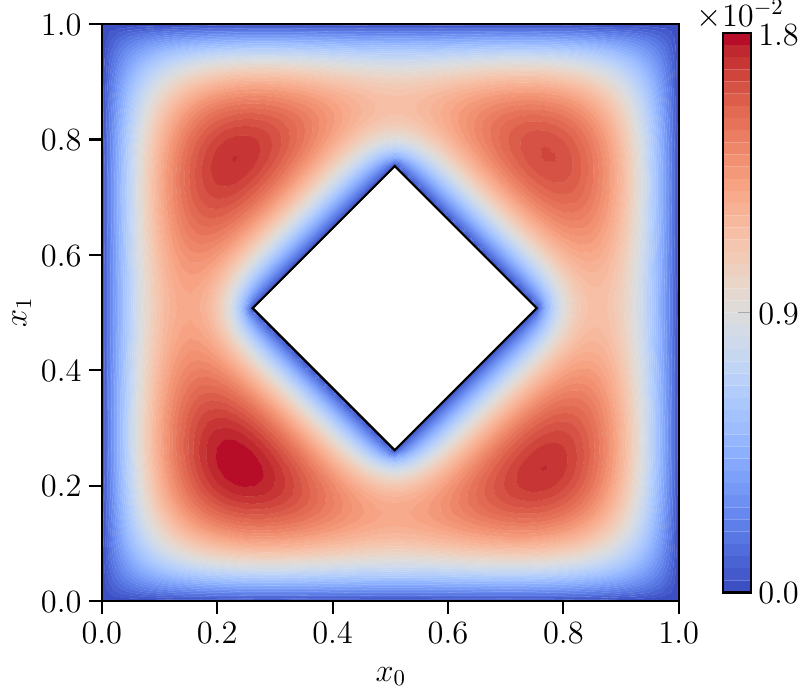}
        \end{minipage}
    }
    \hspace{0.075\textwidth}
    \subfloat[Absolute error~$|u(x_0,x_1) - u_h(x_0,x_1)|$\label{fig:2d-im-figure_b}]{%
        \begin{minipage}[c]{0.43\textwidth}
            \includegraphics[width=\linewidth, keepaspectratio]{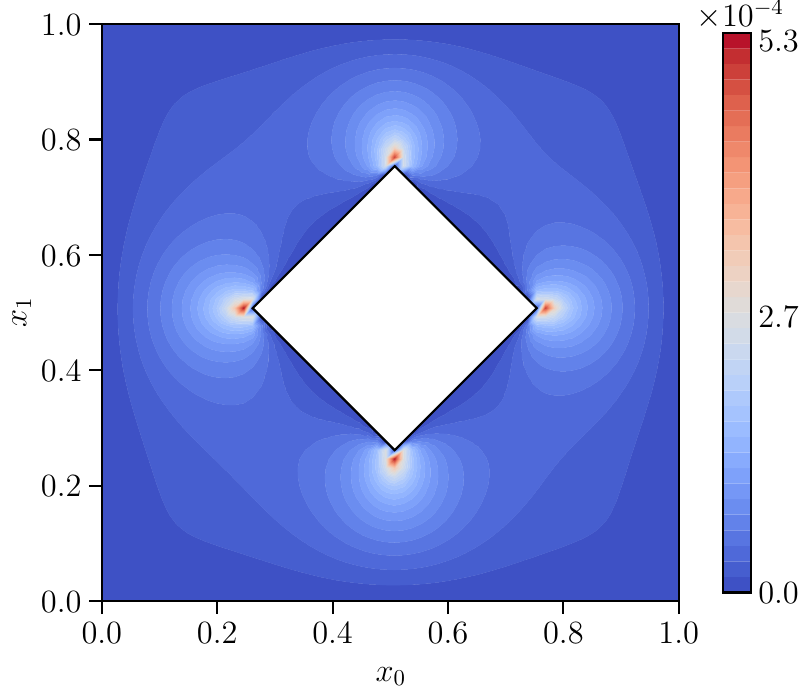}
        \end{minipage}
    } 
    \vspace{0.4cm}
    \subfloat[Number of~$U_3$ and~$\mathit{CNOT}$ gates\label{fig:2d-im-figure_d}]{%
        \begin{minipage}[c]{0.43\textwidth}
            \includegraphics[width=\linewidth, keepaspectratio]{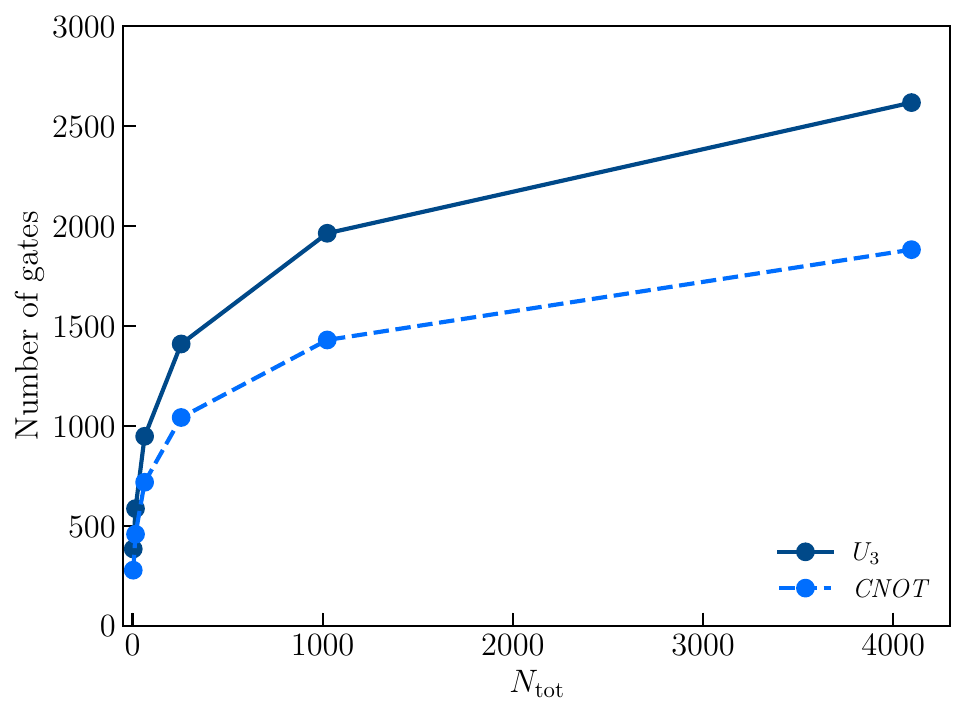}
        \end{minipage}
    }
    \hspace{0.075\textwidth}
    \subfloat[Relative $L^2$-norm error in~$u_h(x_0, x_1)$\label{fig:2d-im-figure_c}]{%
        \begin{minipage}[c]{0.43\textwidth}
            \includegraphics[width=\linewidth, keepaspectratio]{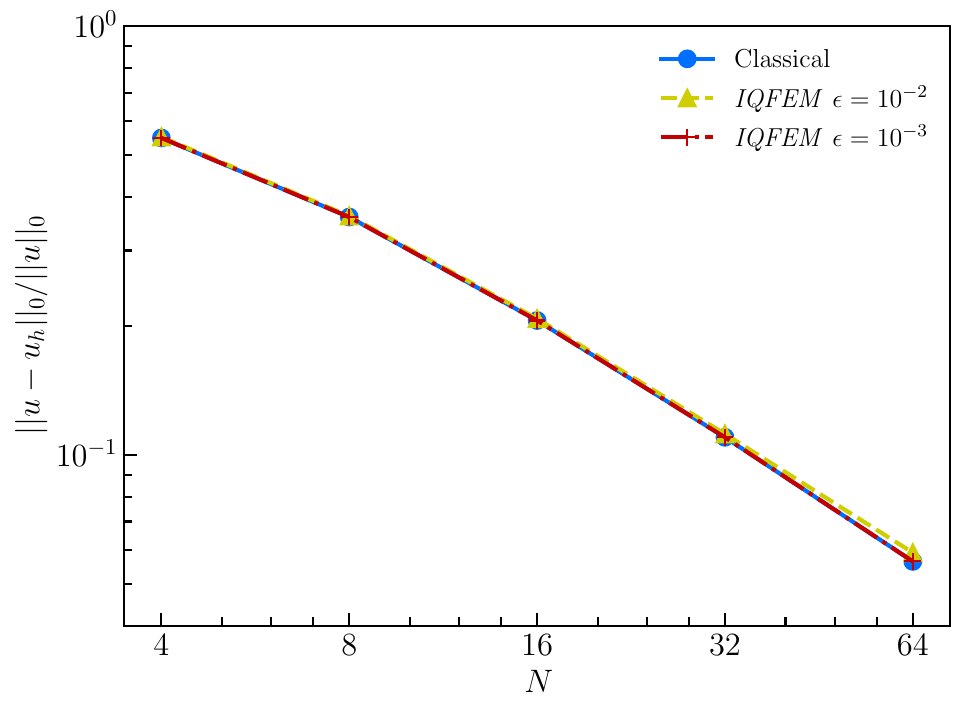}
        \end{minipage}
    }
    \caption{Two-dimensional Poisson problem with immersed Dirichlet boundary. (a) {\em IQFEM} solution~$u_h(x_0, x_1)$ for $2^6 \times 2^6$ grid points. (b) Absolute error~$|u(x_0, x_1) - u_h(x_0, x_1)|$ of the {\em IQFEM} solution compared to a reference fine grid solution. (c) Number of required~$U_3$ and~$\mathit{CNOT}$ gates to block-encode the system matrix. (d) Convergence of the relative $L^2$-norm error of the {\em IQFEM} solution for QSVT approximation errors $\epsilon \in \{10^{-2}, 10^{-3}\}$.}
    \label{fig:2d-im-figure-2}
\end{figure}

The {\em IQFEM} solution on a grid with~$2^6 \times 2^6$ grid points is shown in Figure~\ref{fig:2d-im-figure_a}. For comparison, we compute a reference boundary-fitted classical finite element solution using an unstructured triangular mesh with a characteristic element size of~$ 5\times10^{-4}$. The contours of the pointwise absolute error between {\em IQFEM} and the classical finite element solution are shown in Figure~\ref{fig:2d-im-figure_b}. The errors are localised near the hole boundary and the four sharp corners of the diamond.  These errors arise from the staircase approximation of the hole boundary on the structured Cartesian grid and its inability to represent sharp corners.

The numbers of~$U_3$ and~$\mathit{CNOT}$ gates required to block-encode the immersed system matrix are plotted in Figure~\ref {fig:2d-im-figure_d}. Both gate counts exhibit a polylogarithmic scaling. Their total numbers are slightly higher than those for the preceding simply connected heterogeneous Poisson problem, which can be attributed to the encoding of the level-set and indicator functions.   The convergence of the $L^2$-norm errors in the {\em IQFEM} solution obtained with classical and QSVT-based linear equation solvers and  tolerances~$\epsilon = 10^{-2}$ and~$\epsilon = 10^{-3}$ is shown in Figure~\ref{fig:2d-im-figure_c}. The {\em IQFEM} solution exhibits first-order convergence owing to the staircase approximation of the boundary and corner singularities.

\begin{figure}[!tbp]
   \centering
   \includegraphics[scale=0.43]{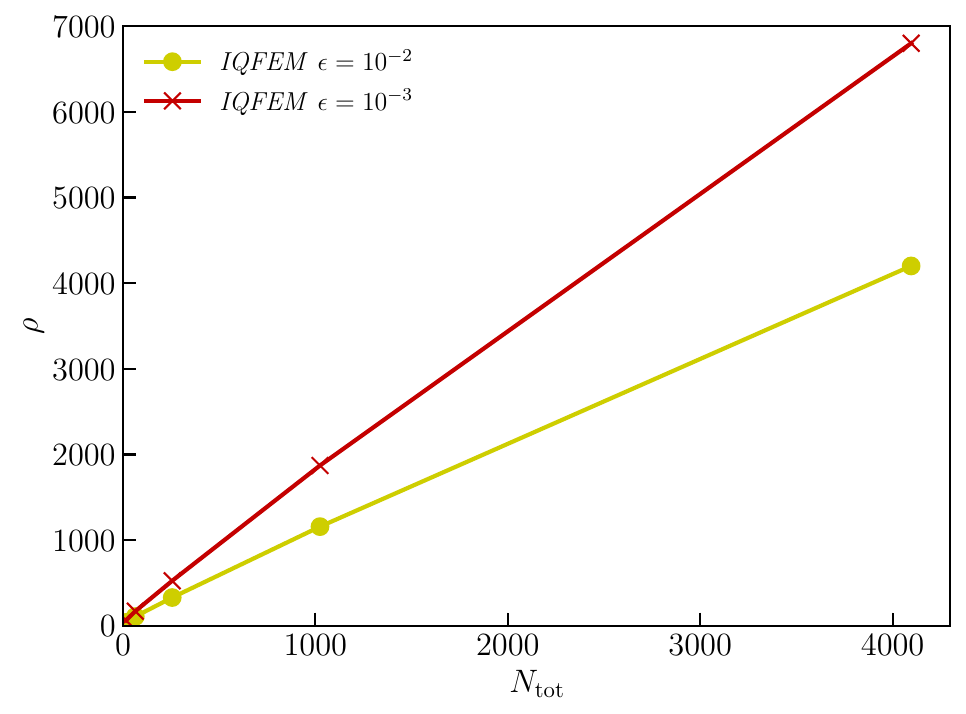}
   \caption{Two-dimensional Poisson problem with immersed Dirichlet boundary. Polynomial degree $\rho$ versus number of grid points for different QSVT approximation errors $\epsilon \in \{10^{-2}, 10^{-3}\}$.\label{fig:2d-im-figure_poly-degree}}
\end{figure}

Figure~\ref{fig:2d-im-figure_poly-degree} shows the Chebyshev polynomial degree~$\rho$ required by the QSVT solver for $N_\text{tot}=N\times N$ grid points. The minimum non-zero eigenvalue of the scaled system matrix $A/ \|A\|_{\max}$ scales as~$N^{-2}$; therefore, grid refinement leads to an approximately linear increase in the required polynomial degree.
%
\subsubsection{Immersed Neumann boundary}
%
Our final example is a multiply connected domain with an immersed six-sided hole. Its geometry within the domain~$(0,  \, 1) \times (0, \, 1)$ is described by the level-set function
\begin{equation}
	\phi (x_0, x_1) = \max \left \{ | x_0 - 1/2| -1/4 , \,   | x_1 - 1/2| -1/4, \, |x_0 + x_1 -1 | -3/8   \right \}  \, ,
\end{equation}
which is constructed from three wedges, each of which is constructed from two half-planes. The diffusion coefficient~$\mu(x_0,x_1)=1$ and the source~$f(x_0,x_1)=1$ are constant.  The external boundary is subject to homogeneous Dirichlet boundary conditions, and the hole boundary to homogeneous Neumann boundary conditions. 

\begin{figure}[!tb]
    \centering
    \subfloat[{\em IQFEM} solution~$u_h(x_0, x_1)$\label{fig:2d-immN-figure_a}]{%
        \begin{minipage}[c]{0.43\textwidth}
            \includegraphics[width=\textwidth, keepaspectratio]{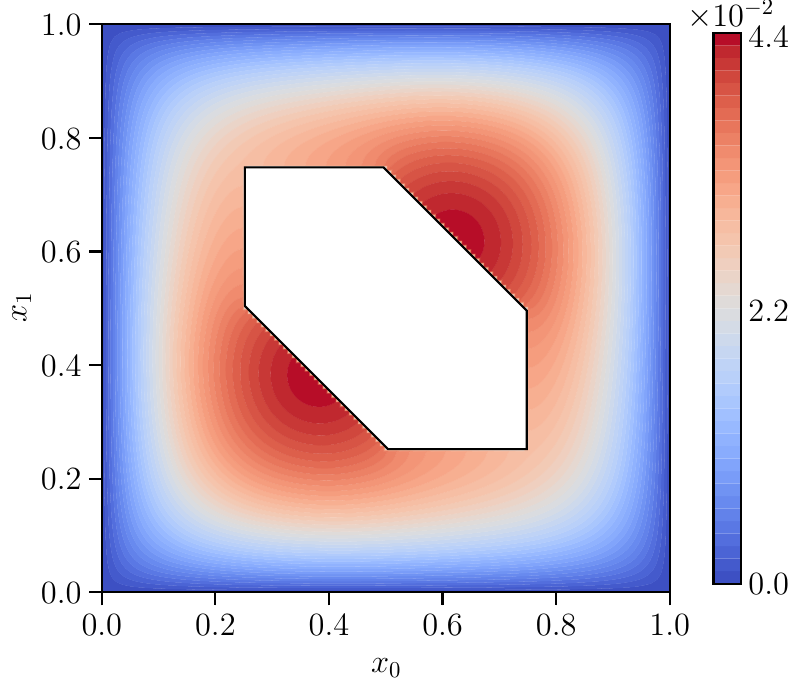}
        \end{minipage}
    }
    \hspace{0.075\textwidth}
    \subfloat[Absolute error~$|u(x_0, x_1) - u_h(x_0, x_1)|$\label{fig:2d-immN-figure_b}]{%
        \begin{minipage}[c]{0.43\textwidth}
            \includegraphics[width=\textwidth, keepaspectratio]{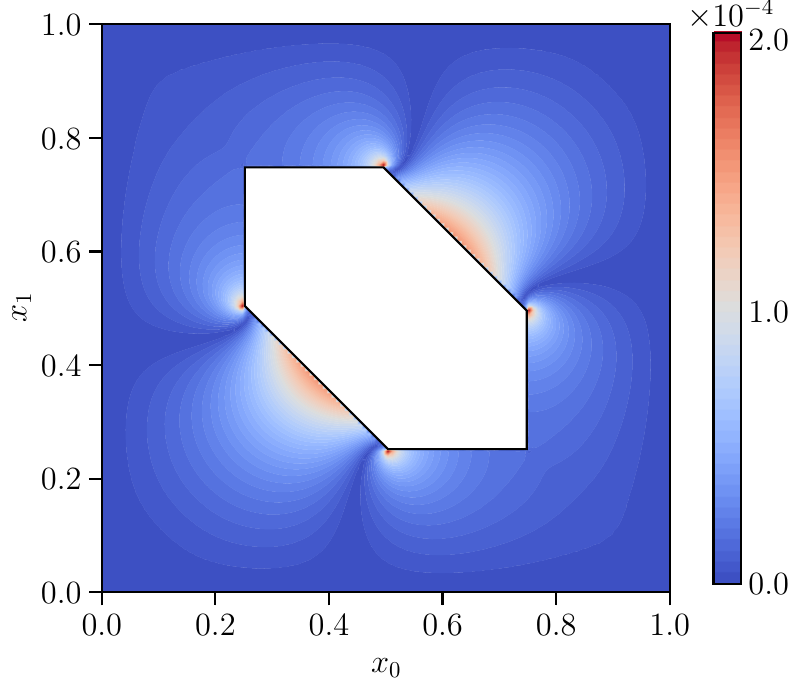}
        \end{minipage}
    }
    \vspace{0.4cm}
    \subfloat[Number of~$U_3$ and~$\mathit{CNOT}$ gates \label{fig:2d-immN-figure_d}]{
        \begin{minipage}[c]{0.43\textwidth}
            \includegraphics[width=\textwidth, keepaspectratio]{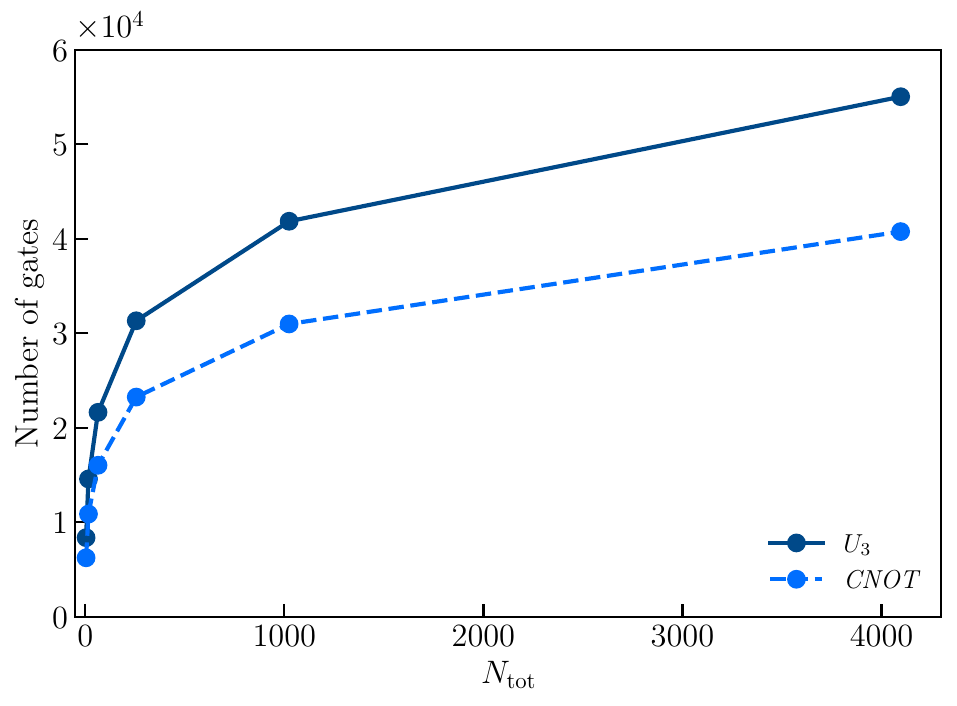}
        \end{minipage}
    }
    \hspace{0.075\textwidth}
    \subfloat[Relative $L^2$-norm error in solution~$u_h(x_0, x_1)$\label{fig:2d-immN-figure_c}]{
        \begin{minipage}[c]{0.42\textwidth}
            \includegraphics[width=\textwidth, keepaspectratio]{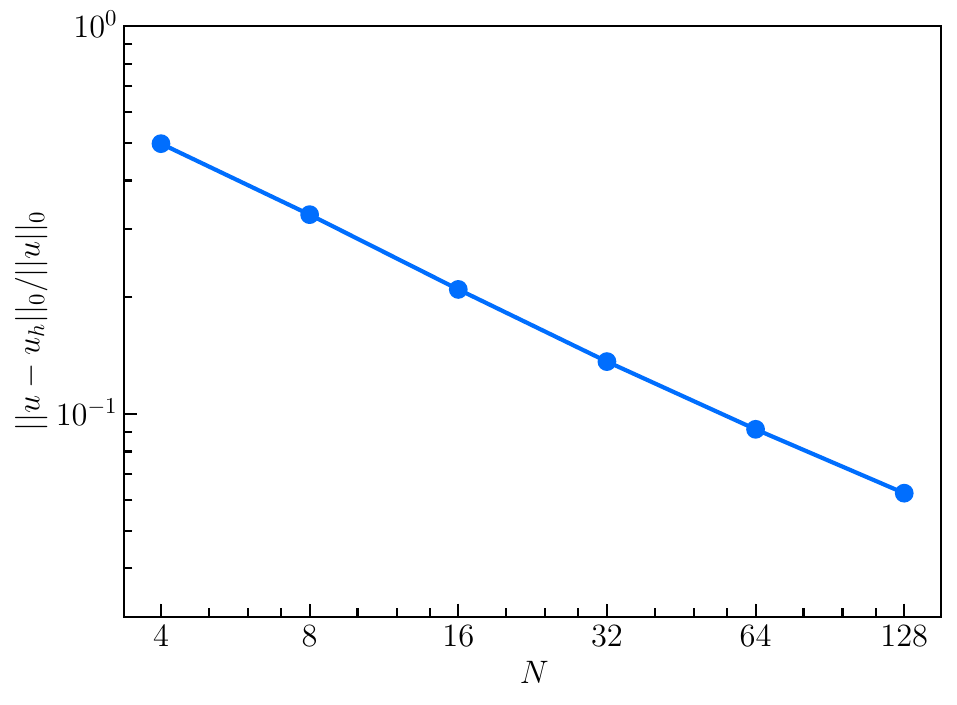}
        \end{minipage}
    }
    \caption{Two-dimensional Poisson problem with immersed Neumann boundary. (a) {\em IQFEM} solution~$u_h(x_0, x_1)$ for $N_\text{tot} = 2^7 \times 2^7$ grid points. (b) Absolute error~$|u(x_0, x_1) - u_h(x_0, x_1)|$  of the {\em IQFEM} solution compared to a reference fine grid FE solution. (c) Number of required~$U_3$ and~$\mathit {CNOT} $ gates to block-encode the system matrix.  (d) Convergence of the relative $L^2$-norm error of the {\em IQFEM} solution. }
    \label{fig:2d-immN-figure}
\end{figure}

The {\em IQFEM} solution on a grid with  $2^7 \times 2^7$ grid points is shown in Figure~\ref{fig:2d-immN-figure_a}. As in the preceding example, we compute a classical reference solution using the boundary-fitted finite element method. The contours of the pointwise absolute error between the {\em IQFEM} and classical finite element solutions are shown in Figure~\ref{fig:2d-immN-figure_b}. The errors near the hole boundary, arising from the staircase approximation, are evident. 

The~$U_3$ and~$\mathit{CNOT}$ gate counts required to block-encode the system matrix are shown in Figure~\ref{fig:2d-immN-figure_d}. Both gate counts exhibit polylogarithmic scaling. The total number of gates is significantly higher than in the previous examples,  primarily due to the use of LCU in computing the immersed Neumann matrix \mbox{$A_N= A_D - D_N(I \circ A)/2d$}, cf. Section~\ref{sec:immersed_sysmat_N}. The block-encodings of the two summands are combined using LCU, which requires controlled versions of both unitaries. Introducing such controls generally leads to a significant increase in the gate counts, although many of them can often be optimised away. Figure~\ref{fig:2d-immN-figure_c} depicts the convergence of the relative $L^2$-norm error, which decreases monotonically at a rate slightly below first order. 

%
\section{Conclusions \label{sec:conclusions}}
%
We have introduced {\em IQFEM} for solving heterogeneous Poisson problems on general domains. The central contribution of {\em IQFEM} is the efficient block encoding of system matrices with a gate count that scales polylogarithmically with the number of grid points. Problems are discretised using structured Cartesian grids, leading to highly structured system matrices whose entries are polynomials in the row indices and, hence, the corresponding grid-point coordinates. We express such system matrices as linear combinations of block-encoded diagonal matrices with polynomial entries and their products with block-encoded shift operators. Their sum is obtained using an LCU construction. The resulting block encoding of the system matrix has a modest number of ancillas for summation and polynomial encoding.  Immersed boundaries are described by integer-valued level-set functions constructed using a solid-geometry approach and evaluated on the fly using quantum arithmetic. The block-encoded immersed system matrix is obtained by pre- and post-multiplying the standard block-encoded system matrix with a level-set-defined diagonal indicator matrix and applying certain corrections. 

We solve the block-encoded linear system of equations using QSVT. The matrix inverse is determined by first approximating $1/x$, after regularisation, by a polynomial and then applying this polynomial via QSVT to the singular values of the system matrix. The required degree of the polynomial increases as the smallest non-zero singular value decreases. In our numerical experiments, approximations using high-degree polynomials, reaching several thousands, did not lead to stability problems owing to the unitarity of the QSVT construction.   The number of invocations of the block-encoded system matrix in QSVT, and hence the total number of elementary gates, is proportional to the polynomial degree. Assuming that the condition number scales as $\kappa\sim N_{\text{tot}}^{2/d}$, the overall gate count scales as $O\left(N_{\text{tot}}^{2/d}\operatorname{polylog}(N_{\text{tot}})\right)$ in terms of the total number of grid points~$N_{\text{tot}}=N^d$. Consequently, without preconditioning, the total number of elementary gates in {\em IQFEM} scales quadratically in 1D, linearly in 2D and sublinearly in 3D. 

The numerical experiments for one- and two-dimensional Poisson problems with spatially varying coefficients and domains with an internal hole confirm both the accuracy and the predicted gate count scaling. Solutions of problems on simply connected domains exhibit the optimal quadratic convergence rate, whereas solutions of immersed boundary problems exhibit suboptimal linear convergence rates. The latter is a consequence of the staircase-like approximation of the immersed boundaries using the integer-valued level-set function. The geometric approximation may be improved by representing the level-set function on a finer subgrid and accounting for the contributions of cut elements. Such an extension would require fixed- or floating-point quantum arithmetic and would increase the qubit and gate counts associated with processing the level-set function and computing the subgrid contributions to the system matrix.

A notable limitation and open challenge for {\em IQFEM} is the need for preconditioning of the system matrix, which is key to achieving better-than-linear gate count scaling for the entire solution procedure. Resource-efficient quantum preconditioning approaches are an active area of research, with some limited work so far~\cite{ shao2018quantum,tong2021fast,golden2022quantum, deiml2024quantum,lapworth2025preconditioned}. Some of these approaches, like the BPX preconditioner in~\cite{deiml2024quantum} or the circulant preconditioner in~\cite{ shao2018quantum}, appear particularly compatible with {\em IQFEM} because of its use of structured grids. Their effectiveness, especially in the presence of immersed boundaries, requires investigation. Finally, to achieve an end-to-end quantum advantage in realistic engineering problems, the readout of the quantum solution must be taken into consideration. Depending on the required output, techniques such as amplitude amplification~\cite{Brassard2002} and classical shadows~\cite{huang2020predicting} may be used to reduce the associated sampling cost.


\appendix
%
\section{Encoding of circulant system matrices \label{app:laplacian}}
%
We consider the quantum encoding of the system matrix~$A$ of a one-dimensional homogeneous periodic Poisson problem discretised with a uniform grid of elements of unit length. The associated circulant tridiagonal matrix~\mbox{$A \in \mathbb R^{N\times N}$} decomposes into
\begin{equation} \label{eq:Laplacedecompose}
    A =  	\begin{pmatrix*}
		2 & -1 &  & & -1  \\
		-1 & 2 & -1 & &   \\
		& -1 & \ddots& \ddots &  \\
		& & \ddots & \ddots & -1    \\
		-1 &  &   &  -1 & 2   
	\end{pmatrix*} 
 =  2I - \tilde S_{+} - \tilde S_{-} \, ,
\end{equation}
where $I$ is the identity matrix and~$\tilde S_{+}$ and~$\tilde S_{-}$ are two cyclic shift operators defined as 
\begin{equation} \label{eq:shiftOperators_cyclic}
	 \tilde S_{+} 
		\begin{pmatrix*}
		v_0  \\
		v_1  \\
		 \vdots    \\
		 v_{N-2} \\
		v_{N-1}    
	\end{pmatrix*} 
	= \begin{pmatrix*} v_{N-1} \\ v_0  \\ \vdots    \\v_{N-3} \\v_{N-2}    \end{pmatrix*}
	\, , \qquad    \tilde S_{-}   
		\begin{pmatrix*}
		v_0  \\
		v_1  \\
		 \vdots    \\
		 v_{N-2} \\
		v_{N-1}    
	\end{pmatrix*}   = 
	\begin{pmatrix*} v_1  \\ v_2  \\  \vdots    \\  v_{N-1}\\ v_0    \end{pmatrix*}  \, . 
\end{equation}

\begin{figure*}[!tbp]
    \centering
    \scalebox{1}{
    	\includegraphics[scale=1]{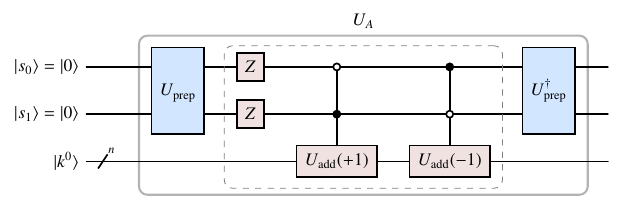} 
    }
    \caption{Block encoding~$U_A$ of the circulant tridiagonal system matrix~$A= 2 I - \tilde S_+ - \tilde S_-$ of a one-dimensional homogeneous periodic problem. The discretisation consists of $N$ cells and grid points, i.e.~$N=2^n$. The two~$U_{\text{add}}$ gates implement the up and down shifts of the basis components. In the output, the relevant components correspond to basis states of the form $\ket{00}\ket{k^0}$. \label{fig:1DlapCirc}}
\end{figure*}

The block-encoding of~$A$ using the LCU technique is depicted in Figure~\ref{fig:1DlapCirc}. The resulting unitary acts on an enlarged space with~$n+2$ qubits, although only the subspace in which the ancillas~$\ket{s_0 s_1}$ are in state~$\ket{00}$ is relevant.
The state preparation unitary~$U_{\text{prep}}$ depends on the prefactors~$(2, \, -1, \, -1)$ in the decomposition~\eqref{eq:Laplacedecompose}. Specifically, it prepares the state~$( 1/\sqrt{2}, \, 1/2, \, 1/2 )$,  which is equal to the square roots of the absolute coefficients normalised by their $L^1$ norm, i.e.~$2 + 1+ 1=4$.  The first shaded gate in Figure~\ref{fig:1DlapCirc} implements the unitary 
\begin{equation}
	U_{\text{prep}} \otimes I =  \begin{pmatrix} \frac{1}{\sqrt 2} I  &  * &  * & * \\[0.5em]  \frac{1}{2}  I &  * &  * & *  \\[0.5em]  \frac{1}{2} I &  * &  * & * \\[0.5em]  0 &  * &  * & * \end{pmatrix} \, .
\end{equation}
Notice that the coefficients in the first column are equal to the scaled positive prefactors and the remaining coefficients are chosen arbitrarily to satisfy unitarity.  The subsequent gates in the dashed box in Figure~\ref{fig:1DlapCirc} implement the select unitary 
\begin{equation}
	U_{\text{select}} =  \begin{pmatrix} I  &   &   &  \\[0.5em]   &  -U_{\text{add}} (+1) &   &   \\[0.5em]  &   &  -U_{\text{add}}(-1) &  \\[0.5em]   &   &   & I \end{pmatrix} \, .
\end{equation}
Here, the unitaries~$U_{\text{add}}(+1)$ and~$U_{\text {add}}(-1)$ implement the cyclic shift operators~${\tilde S}_{+}$ and~${\tilde S}_{-}$, respectively. The negative signs are introduced by the two Pauli Z gates.  The  entire circuit provides the following block encoding for the scaled system matrix~$A/4$,
\begin{equation}
\begingroup
\setlength{\arraycolsep}{3.75pt}
	\begin{pmatrix} \frac{1}{\sqrt 2} I  &   \frac{1}{2}  I   &  \frac{1}{2} I  & 0 \\[0.5em]  *  &  * &  * & *  \\[0.5em]  * &  * &  * & * \\[0.5em]  * &  * &  * & * \end{pmatrix} 
	\begin{pmatrix} I  &   &   &  \\[0.5em]   &  -U_{\text{add}} (+1) &   &   \\[0.5em]  &   &  -U_{\text{add}}(-1) &  \\[0.5em]   &   &   & I \end{pmatrix}
	\begin{pmatrix} \frac{1}{\sqrt 2} I  &  * &  * & * \\[0.5em]  \frac{1}{2}  I &  * &  * & *  \\[0.5em]  \frac{1}{2} I &  * &  * & * \\[0.5em]  0 &  * &  * & * \end{pmatrix}
	= \frac{1}{4} \begin{pmatrix} 2 I - U_{\text{add}}(+1)  - U_{\text{add}}(-1)   & * & *  & * \\ * & * & * & *  \\ * & * & * & * \\ * & * & * & * \end{pmatrix} \, .
\endgroup
\end{equation}
Hence, for the input state~$\ket  {00} \ket {v}$, the component of the output in the subspace with the top two qubits in state~$\ket {00}$ is~$(A/4) \ket v$. The desired output~$A \ket v $ is obtained, up to normalisation, by postselecting this subspace. The corresponding success probability is $\| (A/4) \ket {v}\|^2$. See \cite{camps2024explicit,sunderhauf2024block} for further details.

%
\section{QFT-based quantum arithmetic \label{app:arithmetic}}
%
There are several approaches for implementing integer and fixed-point arithmetic for binary-encoded numbers on a quantum computer. Such an encoding of numbers is also referred to as basis encoding.  Numbers are quantum encoded by setting each qubit's state to either~$0$ or~$1$ according to their binary representation. We review the QFT-based approach to quantum arithmetic for integers and briefly comment on its application to fixed-point numbers~\cite{draper2000addition,ruiz2017quantum,wang2025comprehensive}.

%
\subsection{Addition \label{app:addition}}
%
Two different addition operations are used throughout the paper, depending on whether both summands are quantum or only one is quantum. We refer to the former as quantum-quantum addition and the latter as quantum-classical addition.
%
\subsubsection{Quantum-classical addition \label{app:qc-addition}}
\begin{figure}[!tb]
    \centering
    \includegraphics[scale=1]{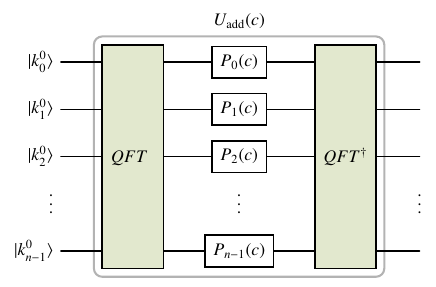} 
    \caption{Unitary $U_{\text{add}}(c)$ for adding the classical integer~$c$ and the quantum basis-encoded integer~$k$. \label{fig:qft-shifter}}
\end{figure}

%
Let $c,k \in \{0, \, 1, \, \dotsc, \, 2^n-1  \}$, where~$c$ is a given classical integer and~$k$ is the quantum basis-encoded integer 
\begin{equation}
k = \sum_{j=0}^{n-1} k_j 2^{\,n-1-j} \, . 
\end{equation}
We consider the construction of~$U_{\text{add}}(c)$ implementing the modular addition
\begin{equation}
    U_{\text{add}}(c) \colon \ket{k} \;\mapsto\; \ket{(k+c)\bmod 2^n} \, .
\end{equation}
It is expedient to work in the Fourier basis to implement this unitary. The Fourier transformed state $\mathit{QFT}\ket{k}$ corresponds to the $k$-th column of the unitary Fourier matrix.  Adding $c$ amounts to shifting from the $k$-th to the $(k+c)$-th column, and the columns of the Fourier matrix differ only by phase factors. The key observation is that the column indexed by  $k$ admits the tensor-product representation
\begin{equation}
\mathit{QFT} \ket{k}
= \frac{1}{\sqrt{2^n}}\bigotimes_{j=0}^{n-1} 
\begin{pmatrix}
1 \\ e^{2\pi i k/2^{j+1}}
\end{pmatrix} \, ,
\end{equation} 
where the two-dimensional vector represents the state of qubit~$j$. Replacing the integer $k$ by $k+c$ yields the column
\begin{equation}
\mathit{QFT} \ket{k+c}
= \frac{1}{\sqrt{2^n}}\bigotimes_{j=0}^{n-1}
\begin{pmatrix}
1 \\ e^{2\pi i (k+c)/2^{j+1}}
\end{pmatrix} \, .
\end{equation}
Noting that $e^{2\pi i (k+c)/2^{j+1}} = e^{2\pi i k/2^{j+1}} e^{2\pi i c/2^{j+1}}$, the mapping of the state of the qubit~$j$ can be realised via the single-qubit phase gate
\begin{equation}\label{eq:phase-gate}
P_j(c) =
\begin{pmatrix}
1 & 0 \\
0 & e^{2\pi i c/2^{j+1}}
\end{pmatrix}\, .
\end{equation}
The unitary for modular addition can thus be decomposed into a layer of phase rotations sandwiched between the $\mathit{QFT}$ and its inverse, 
\begin{equation}
U_{\text{add}}(c) = \mathit{QFT}^\dagger
\left(
\bigotimes_{j=0}^{n-1} P_j(c)\right)\mathit{QFT}\, .
\end{equation}
See Figure~\ref{fig:qft-shifter} for the quantum implementation of this unitary. 

%
\subsubsection{Quantum-quantum addition \label{app:qq-addition}}
%
Let $k^0,k^1 \in \{0, \, 1, \, \dotsc, \, 2^n-1  \}$ be two quantum basis-encoded integers 
\begin{equation}
k^0 = \sum_{j=0}^{n-1} k^0_j 2^{\,n-1-j}\, ,
\qquad
k^1 = \sum_{j=0}^{n-1} k^1_j 2^{\,n-1-j}\, .
\end{equation}
We consider the construction of~$U_{\text{qadd}}$ for the addition of these two integers. To implement non-modular addition, the first register is augmented by one qubit and the composite quantum state is initialised as $\ket{0} \ket{k^0} \ket{k^1}$. The additional qubit stores the carry, defined as 
\begin{equation}
\mathrm{carry}(k^0,k^1) = \left\lfloor \frac{k^0 + k^1}{2^n} \right\rfloor \, .
\end{equation}
The unitary~$U_{\text{qadd}}$ implements the in-place addition 
\begin{equation}
U_{\text{qadd}} \colon \ket{0} \ket{k^0} \ket{k^1}
\;\mapsto\;
\ket{\mathrm{carry}(k^0,k^1)}
\, \ket{(k^0 + k^1)\bmod 2^n}
\, \ket{k^1}\, .
\end{equation}

\begin{figure}[ht]
    \centering
    \includegraphics[scale=1]{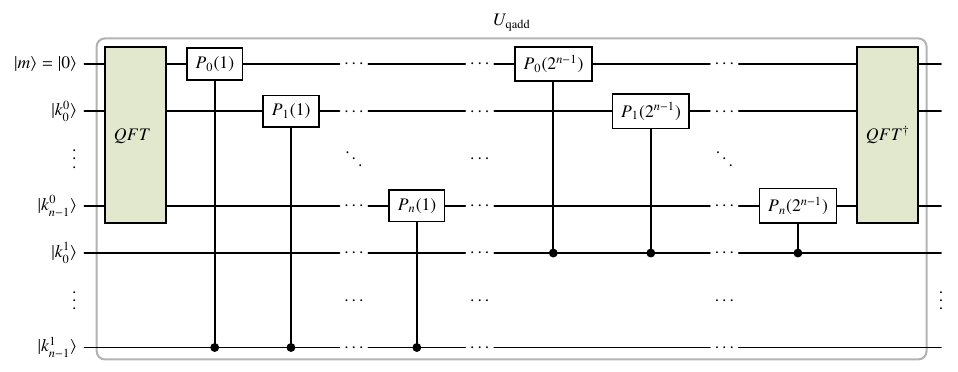} 
    \caption{Unitary $U_{\text{qadd}}$ for in-place addition of two quantum basis-encoded integers~$k^0$ and~$k^1$. \label{fig:full_adder_quantum}}
\end{figure}

Again, it is convenient to work in the Fourier basis. Applying $\mathit{QFT}$ to the register~$\ket{0}\ket{k^0}$ maps~$k^0$ to the corresponding column of the QFT matrix. Adding $k^1$ amounts to a column shift, which factorises into phase rotations on each qubit. Each bit $k^1_j$ in the binary representation of~$k^1$ induces a controlled phase rotation on qubit $r$ with angle
\begin{equation}
\theta_{rj} = \frac{2\pi \, 2^{\,n-1-j}}{2^{r+1}}\, .
\end{equation}
After applying the controlled rotations and~$\mathit{QFT}^\dagger$, the register~$\ket{0}\ket{k^0}$ contains the full $(n+1)$-bit sum, while the register~$\ket{k^1}$ remains unchanged. See Figure~\ref{fig:full_adder_quantum} for the quantum implementation of this approach. 

%
\subsection{Absolute value \label{app:absolute}}
%
We consider the construction of the unitary~$U_{\text{abs}}$  implementing the mapping $k \mapsto \vert k-c \vert$, where~$k$ is a quantum basis-encoded integer and~$c$ is a given classical integer. The $n$-qubit register encoding~$\ket{k}$ is augmented by one qubit as~$\ket{0}\ket{k}$.  The additional qubit stores the borrow bit. 

Subtraction is implemented via the introduced quantum-classical addition operator
\begin{equation}
U_{\text{add}}(-c) \colon \ket{0}\ket{k} \;\mapsto\; \ket{m} \ket{d} \, ,
\end{equation}
where $d = k-c \pmod{2^n}$, and $m=1$ if $k<c$ and $m=0$ otherwise.  Hence, 
\begin{equation}
k-c = d - m2^n \, .
\end{equation}
Therefore, the absolute value is given by
\begin{equation}
\vert k-c \vert =
\begin{cases}
d\, , & m=0\, ,\\
2^n - d\, , & m=1 \, .
\end{cases}
\end{equation}
The mapping $d \mapsto 2^n - d$ is applied conditionally on $m=1$ using a bitwise inversion followed by the incrementation $U_{\text{add}}(+1)$. 

\begin{figure}[ht]
    \centering
    \includegraphics[scale=1]{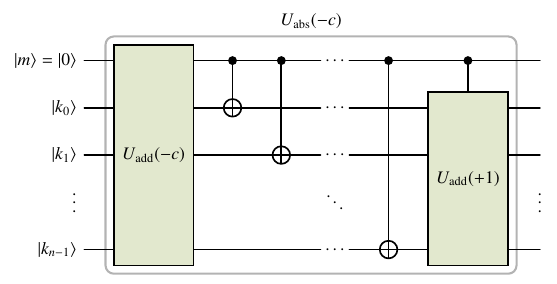} 
    \caption{Unitary $U_{\text{abs}}(-c)$ for computing the absolute value~$|k-c|$ between the quantum basis-encoded integer~$k$ and the classical integer~$c$.   \label{fig:qft-abs}}
\end{figure}

The overall unitary~$U_{\text{abs}}(-c)$ thus implements 
\begin{equation}
	U_{\text{abs}}(-c) \colon \ket{0} \ket{k} \;\mapsto\; \ket{m} \ket{\,\vert k-c \vert\,} \, ;
\end{equation}
see Figure~\ref{fig:qft-abs} for the corresponding quantum circuit.

%
\subsection{Fixed-point numbers \label{app:floating_point}}
%
The unitaries for arithmetic primitives introduced so far directly apply to non-integer fixed-point numbers with prescribed integer and fractional bits. The distinction between integer and fixed-point numbers is purely interpretational and does not affect unitary operations.

Consider the $n$-qubit registers~$\ket{k} = \ket{k_0 k_1 \dots k_{n-1}}$ and~$\ket{l} = \ket{l_0 l_1 \dots l_{n-1}}$. While previously interpreted as an integer, the register may instead be interpreted as a fixed-point number.  For a chosen number of fractional bits $f$, their values in decimal are given by
\begin{equation}
x = 2^{-f} \sum_{j=0}^{n-1} k_j 2^{n-1-j} \, , \quad  y = 2^{-f} \sum_{j=0}^{n-1} l_j 2^{n-1-j} \, . 
\end{equation}
Under this interpretation, the QFT-based arithmetic circuits act on fixed-point numbers without modification. Since $x + y = 2^{-f}(k + l)$, modular addition on the underlying integer encodings preserves the numerical relation, provided the binary points are aligned.

As an illustration, consider adding $k = \ket{01.101}$ (in decimal, $1.625$) and $l = \ket{00.011}$ (in decimal, $0.375$). The quantum-quantum adder yields the state $\ket{10.000}$, which corresponds to $2.0$. This demonstrates that the adder correctly propagates carries independently of whether the bit string is interpreted as an integer or a scaled fractional quantity.

%
\section{Products of block-encoded matrices \label{app:product_encoding}}
%
The product of two block-encoded matrices~$B,C\in\mathbb{R}^{N\times N}$ can be computed straightforwardly when separate ancilla registers are used to encode~$B$ and~$C$. For simplicity, we assume that the block encodings have unit normalisation factors and that each matrix is block-encoded using a single ancilla qubit, so that $U_B,U_C\in\mathbb{R}^{2N\times 2N}$. Specifically,
\begin{equation}
    U_B \ket 0 \ket{v}
    =
    \begin{pmatrix}
        B & * \\
        * & *
    \end{pmatrix}
    \begin{pmatrix}
        \ket{v} \\
        0
    \end{pmatrix}
    =
    \begin{pmatrix}
        B\ket{v} \\
        *
    \end{pmatrix} \, ,
    \qquad
    U_C \ket 0 \ket{v} 
    =
    \begin{pmatrix}
        C & * \\
        * & *
    \end{pmatrix}
    \begin{pmatrix}
        \ket{v} \\
        0
    \end{pmatrix}
    =
    \begin{pmatrix}
        C\ket{v} \\
        *
    \end{pmatrix} \, .
\end{equation}
For each block encoding, the ancilla is initially in the state~$\ket{0}$, and only the output component with the ancilla in the state~$\ket{0}$ is relevant.  

\begin{figure}
\tikzset{operator/.append style={fill opacity=0.875}}
\centering
\begin{quantikz}
    [row sep={.8cm,between origins}, column sep=0.35cm]
    \lstick{$ \ket{0}$} & \qw & \qw& \gate[3, label style={yshift=2mm}, style={fill=csmlPink!40} ]{U_C} & \qw \\
    \lstick{$\ket{0}$} & \qw& \gate[2, style={fill=csmlBrightBlue!20}]{U_B} &  \linethrough& \qw \\
    \lstick{$\ket{v}$} & \qwbundle{n}  &  & & \qw
\end{quantikz}
\caption{Product of two block-encoded matrices $B$ and~$C$. \label{fig:block_encoding_product}} 
\end{figure}
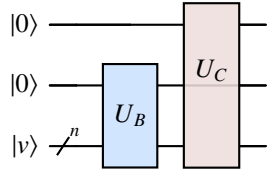
To compute~$CB\ket{v}$, we consider the circuit in Figure~\ref{fig:block_encoding_product}.  The first and second qubits from the top are the ancillas used to encode~$C$ and~$B$, respectively. The
block encodings~$U_B$ and~$U_C$ are applied successively to the~$\ket v$ register using their respective ancillas. The action of the circuit on the state $\ket{0}\ket{0}\ket{v}$ is
\begin{equation}
    \begin{pmatrix}
        C & 0 & * & 0 \\
        0 & C & 0 & * \\
        * & 0 & * & 0 \\
        0 & * & 0 & *
    \end{pmatrix}
        \begin{pmatrix}
        B & * & 0 & 0 \\
        * & * & 0 & 0 \\
        0 & 0 & B & * \\
        0 & 0 & * & *
    \end{pmatrix}
    \begin{pmatrix}
        \ket{v} \\
        0 \\
        0 \\
        0
    \end{pmatrix}
    =
    \begin{pmatrix}
        CB\ket{v} \\
        * \\
        * \\
        *
    \end{pmatrix} \, . 
\end{equation}
The output component with both ancillas in the state~$\ket 0 \ket 0$ yields the desired matrix--vector product~$CB\ket{v}$. Thus, the circuit provides a block encoding of the product~$CB$.

\bibliographystyle{elsarticle-num-names}
\bibliography{iqFEM}

\end{document}